\documentclass[english,aps, schowpacs, prl, reprint, longbibliography]{revtex4-1}
\usepackage[T1]{fontenc}
\usepackage[latin9]{inputenc}
\usepackage{color}
\usepackage{babel}
\usepackage{amsmath}
\usepackage{amssymb}
\usepackage{graphicx}
\usepackage{esint}
\usepackage[unicode=true,
 bookmarks=true,bookmarksnumbered=false,bookmarksopen=false,
 breaklinks=true,pdfborder={0 0 0},backref=false,colorlinks=true]
 {hyperref}
\hypersetup{
 pdfauthor={Jakob Löber}}
\usepackage{breakurl}

\makeatletter
 \@ifundefined{textcolor}{}
 {%
   \definecolor{BLACK}{gray}{0}
   \definecolor{WHITE}{gray}{1}
   \definecolor{RED}{rgb}{1,0,0}
   \definecolor{GREEN}{rgb}{0,1,0}
   \definecolor{BLUE}{rgb}{0,0,1}
   \definecolor{CYAN}{cmyk}{1,0,0,0}
   \definecolor{MAGENTA}{cmyk}{0,1,0,0}
   \definecolor{YELLOW}{cmyk}{0,0,1,0}
 }

\usepackage{tabularx}

\makeatother

\begin{document}

\title{Linear structures in nonlinear optimal control}

\author{Jakob Löber}

\email{jakob@physik.tu-berlin.de}

\selectlanguage{english}%

\affiliation{Institut für Theoretische Physik, EW 7-1, Hardenbergstraße 36, Technische
Universität Berlin, 10623 Berlin}

\date{\today}
\begin{abstract}
We investigate optimal control of dynamical systems which are affine,
i.e., linear in control, but nonlinear in state. The control task
is to enforce the system state to follow a prescribed desired trajectory
as closely as possible, a task also known as optimal trajectory tracking.
To obtain well-behaved solutions to optimal control, a regularization
term with coefficient $\varepsilon$ must be included in the cost
functional. Assuming $\varepsilon$ to be small, we reinterpret affine
optimal control problems as singularly perturbed differential equations.
Performing a singular perturbation expansion, approximations for the
optimal tracking of arbitrary desired trajectories are derived. For
$\varepsilon=0$, the state trajectory may become discontinuous, and
the control may diverge. On the other hand, the analytical treatment
becomes exact. We identify the conditions leading to linear evolution
equations. These result in exact analytical solutions for an entire
class of nonlinear trajectory tracking problems. The class comprises,
among others, mechanical control systems in one spatial dimension
and the FitzHugh-Nagumo model with a control acting on the activator.
\end{abstract}

\pacs{02.30.Yy, 05.45.-a, 07.05.Dz, 45.80.+r}

\maketitle

\paragraph{Introduction.}

In principle, a controlled dynamical systems can be considerably simpler
than the corresponding uncontrolled system. Consider Newton's equation
of motion (EOM), cast in form of a two-dimensional dynamical system,
for a point particle with unit mass, position $x$, and velocity $y$.
Under the influence of an external force $R$ and a control force
$bu$, the EOM read 
\begin{align}
\dot{x} & =y, & \dot{y} & =R\left(x,y\right)+b\left(x,y\right)u.\label{eq:xyState}
\end{align}
While no simple analytical expression for the solution exists for
the uncontrolled system with $u\equiv0$, assuming that $u$ is a
feedback control signal renders a linearization of Eq. (\ref{eq:xyState})
possible. Indeed, for $b\left(x,y\right)\neq0$, we may introduce
a new control signal $v$ by $u=\frac{1}{b\left(x,y\right)}\left(v-R\left(x,y\right)\right)$
such that the new controlled system is linear, $\dot{x}=y,\,\dot{y}=v$.
This technique is called feedback linearization and is applicable,
in a more sophisticated fashion involving an additional state transformation,
to a huge class of dynamical systems \cite{khalil2002nonlinear}.
The resulting linear structure facilitates the application of a multitude
of solution and analysis techniques, which are not available in the
nonlinear case \cite{chen1995linear}. Contrary to the more familiar
approximate linearizations applied e.g. to study of the linear stability
of attractors, feedback linearization is an example of an exact linearization.
The question emerges if exact linearizations are exclusive to systems
with feedback control.

In this letter, we demonstrate the possibility of linear structures
in optimal control, which may or may not be a feedback control. As
a result, the controlled state as well as the control signal are given
as the solution to linear differential and algebraic equations. The
finding enables the exact analytical solution of an entire class of
nonlinear optimally controlled dynamical systems, including but not
limited to all models of the form Eq. (\ref{eq:xyState}). While stabilizing
feedback control of unstable attractors received much attention by
the physics community \cite{bechhoefer2005feedback}, especially in
the context of chaos control \cite{ott1990controlling}, optimality
of these methods is rarely investigated. Numerical solutions of optimal
control problems are computationally expensive, which often prevents
real time applications. Particularly optimal feedback control suffers
from the 'curse of dimensionality' \cite{bellman1957dynamic}. On
the other hand, analytical methods are largely restricted to linear
systems which lack e.g. limit cycles and chaos. The approach presented
here opens up a way to circumvent such problems.

Trajectory tracking aims at enforcing, via a vector of control signals
$\boldsymbol{u}\in\mathbb{R}^{p}$, a system state $\boldsymbol{x}\in\mathbb{R}^{n}$
to follow a prescribed desired trajectory $\boldsymbol{x}_{d}\in\mathbb{R}^{n}$
as closely as possible within a time interval $t_{0}\leq t\leq t_{1}$.
The distance between $\boldsymbol{x}$ and $\boldsymbol{x}_{d}$ in
function space can be measured by the functional $\mathcal{J}=\frac{1}{2}\intop_{t_{0}}^{t_{1}}dt\left(\boldsymbol{x}-\boldsymbol{x}_{d}\right)^{2}$.
$\mathcal{J}$ vanishes if and only if $\boldsymbol{x}=\boldsymbol{x}_{d}$,
in which case we call $\boldsymbol{x}_{d}$ an exactly realizable
desired trajectory \cite{lober2016exactly,loeber2015optimal}. If
$\boldsymbol{x}_{d}$ is not exactly realizable, we can formulate
trajectory tracking as an optimization problem. The control task is
to minimize the quadratic cost functional 
\begin{align}
\mathcal{J} & =\frac{1}{2}\int_{t_{0}}^{t_{1}}dt\left(\boldsymbol{x}-\boldsymbol{x}_{d}\right)^{T}\boldsymbol{\mathcal{S}}\left(\boldsymbol{x}-\boldsymbol{x}_{d}\right)+\frac{\varepsilon^{2}}{2}\int_{t_{0}}^{t_{1}}dt\boldsymbol{u}^{2},\label{eq:OptimalTrajectoryTrackingFunctionalGeneral}
\end{align}
subject to the dynamic constraints that $\boldsymbol{x}$ evolves
according to the controlled dynamical system 
\begin{align}
\boldsymbol{\dot{x}} & =\boldsymbol{R}\left(\boldsymbol{x}\right)+\boldsymbol{\mathcal{B}}\left(\boldsymbol{x}\right)\boldsymbol{u}, & \boldsymbol{x}\left(t_{0}\right) & =\boldsymbol{x}_{0}, & \boldsymbol{x}\left(t_{1}\right) & =\boldsymbol{x}_{1}.\label{eq:ControlledState}
\end{align}
Here, $\boldsymbol{\mathcal{S}}=\boldsymbol{\mathcal{S}}^{T}$ is
a constant symmetric positive definite $n\times n$ matrix of weights.
While the nonlinearity $\boldsymbol{R}\left(\boldsymbol{x}\right)\in\mathbb{R}^{n}$
is known from uncontrolled dynamical systems, the $n\times p$ input
matrix $\boldsymbol{\mathcal{B}}\left(\boldsymbol{x}\right)$ is exclusive
to controlled systems. We assume that the rank of the matrix $\boldsymbol{\mathcal{B}}\left(\boldsymbol{x}\right)$
equals the number $p$ of independent components of $\boldsymbol{u}$
for all $\boldsymbol{x}$.

The term involving $\varepsilon\geq0$ in Eq. (\ref{eq:OptimalTrajectoryTrackingFunctionalGeneral})
favors controls with small amplitude and serves as a regularization
term. The idea is to use $1\gg\varepsilon\geq0$ as the small parameter
for a perturbation expansion. Note that $\varepsilon$ has its sole
origin in the formulation of the control problem. Assuming it to be
small does not involve simplifying assumptions about the system dynamics.
The dynamics is taken into account without approximations in the subsequent
perturbation expansion. Furthermore, among all optimal controls, the
unregularized ($\varepsilon=0$) one brings the controlled state closest
to the desired trajectory $\boldsymbol{x}_{d}$. For a given dynamical
system, the unregularized optimal control solution can be seen as
the \textit{limit of realizability} of a prescribed desired trajectory
$\boldsymbol{x}_{d}$.

Following a standard procedure \cite{pontryagin1987mathematical,bryson1969applied},
the constrained optimization problem Eqs. (\ref{eq:OptimalTrajectoryTrackingFunctionalGeneral})
and (\ref{eq:ControlledState}) is converted to an unconstrained optimization
problem by introducing the vector of Lagrange multipliers $\boldsymbol{\lambda}\in\mathbb{R}^{n}$,
also called the co-state or adjoint state. Minimizing the constrained
functional with respect to $\boldsymbol{x},\,\boldsymbol{\lambda}$,
and $\boldsymbol{u}$ yields the necessary optimality conditions,
\begin{align}
\boldsymbol{0} & =\varepsilon^{2}\boldsymbol{u}+\boldsymbol{\mathcal{B}}^{T}\left(\boldsymbol{x}\right)\boldsymbol{\lambda},\label{eq:StationarityCondition}\\
\boldsymbol{\dot{x}} & =\boldsymbol{R}\left(\boldsymbol{x}\right)+\boldsymbol{\mathcal{B}}\left(\boldsymbol{x}\right)\boldsymbol{u},\label{eq:StateEq}\\
-\boldsymbol{\dot{\lambda}} & =\left(\nabla\boldsymbol{R}^{T}\left(\boldsymbol{x}\right)+\left(\nabla\boldsymbol{\mathcal{B}}\left(\boldsymbol{x}\right)\boldsymbol{u}\right)^{T}\right)\boldsymbol{\lambda}+\boldsymbol{\mathcal{S}}\left(\boldsymbol{x}-\boldsymbol{x}_{d}\right),\label{eq:AdjointEq}
\end{align}
with $n\times n$ matrix $\left(\nabla\boldsymbol{\mathcal{B}}\left(\boldsymbol{x}\right)\boldsymbol{u}\right)_{ij}=\sum_{k}\partial_{j}\mathcal{B}_{ik}\left(\boldsymbol{x}\right)u_{k}$,
and Jacobi matrix $\left(\nabla\boldsymbol{R}\right)_{ij}=\partial_{j}R_{i}$
of $\boldsymbol{R}$, and $\partial_{j}=\frac{\partial}{\partial x_{j}}$.
The dynamics of an optimal control system takes place in the combined
space of state $\boldsymbol{x}$ and co-state $\boldsymbol{\lambda}$
with dimension $2n$. Starting from Eqs. (\ref{eq:StationarityCondition})-(\ref{eq:AdjointEq}),
the usual procedure is to eliminate $\boldsymbol{u}$ and solve the
system of $2n$ coupled ordinary differential equations (ODEs) for
$\boldsymbol{\lambda}$ and $\boldsymbol{x}$.

To take advantage of the small parameter $\varepsilon$, we proceed
differently. We rearrange the necessary optimality conditions such
that $\varepsilon$ multiplies the highest order derivative of the
system. The rearrangement admits an interpretation of a singular optimal
control problem as a singularly perturbed system of ODEs. Setting
$\varepsilon=0$ yields the \textit{outer equations} but changes the
differential order of the system. Remarkably, the outer equations
may become linear even if the original system is nonlinear. We collect
the conditions leading to linear outer equations under the name \textit{linearizing
assumption.}

Because $\varepsilon=0$ decreases the differential order, not all
$2n$ initial and terminal conditions can be satisfied. Consequently,
the outer solutions are not uniformly valid over the entire time domain.
The non-uniformity manifests itself in initial and terminal boundary
layers of width $\varepsilon$ for certain components of the state
$\boldsymbol{x}$. The boundary layers are resolved by the left and
right \textit{inner equations} valid near to the beginning and the
end of the time interval, respectively. For $\varepsilon>0$, inner
and outer solutions can be composed to an approximate \textit{composite
solution} valid over the entire time domain. See \cite{bender1999advanced}
for analytical methods based on singular perturbations. The control
signal is given in terms of the composite solution. The control exhibits
a maximum amplitude at the positions of the boundary layers. In general,
the inner equations are nonlinear even if the linearizing assumptions
holds. However, for $\varepsilon\rightarrow0$, the boundary layers
degenerate into discontinuous jumps located at the beginning and end
of the time interval. The jump heights are independent of any details
of the inner equations, and the exact solution for $\varepsilon=0$
is governed exclusively by the outer equations. Consequently, the
linearizing assumption together with $\varepsilon=0$ renders the
nonlinear optimal control system linear. On the downside, the optimally
controlled state becomes discontinuous at the time domain boundaries,
and the control diverges. It is in this sense that we are able to
speak about linear structures in nonlinear optimal control.

\paragraph{Rearranging the necessary optimality conditions.}

The rearrangement is based on two complementary $n\times n$ projection
matrices
\begin{align}
\boldsymbol{\mathcal{P}}\left(\boldsymbol{x}\right) & =\boldsymbol{\Omega}\left(\boldsymbol{x}\right)\boldsymbol{\mathcal{S}}, & \boldsymbol{\mathcal{Q}}\left(\boldsymbol{x}\right) & =\mathbf{1}-\boldsymbol{\mathcal{P}}\left(\boldsymbol{x}\right),
\end{align}
with symmetric $n\times n$ matrix 
\begin{align}
\boldsymbol{\Omega}\left(\boldsymbol{x}\right) & =\boldsymbol{\mathcal{B}}\left(\boldsymbol{x}\right)\left(\boldsymbol{\mathcal{B}}^{T}\left(\boldsymbol{x}\right)\boldsymbol{\mathcal{S}}\boldsymbol{\mathcal{B}}\left(\boldsymbol{x}\right)\right)^{-1}\boldsymbol{\mathcal{B}}^{T}\left(\boldsymbol{x}\right).\label{eq:DefOmegaGamma}
\end{align}
We drop the dependence on the state $\boldsymbol{x}$ in the following.
It is understood that $\boldsymbol{R}$ and all matrices, except $\boldsymbol{\mathcal{S}}$,
may depend on $\boldsymbol{x}$. The ranks of the projectors $\boldsymbol{\mathcal{P}}$
and $\boldsymbol{\mathcal{Q}}$ are $p$ and $n-p$, respectively,
such that $\boldsymbol{\mathcal{P}}\boldsymbol{x}\in\mathbb{R}^{n}$
has only $p$ linearly independent components. The projectors satisfy
idempotence, $\boldsymbol{\mathcal{P}}^{2}=\boldsymbol{\mathcal{P}}$
and $\boldsymbol{\mathcal{Q}}^{2}=\boldsymbol{\mathcal{Q}}$, and
$\boldsymbol{\mathcal{P}}\boldsymbol{\mathcal{B}}=\boldsymbol{\mathcal{B}}$,
$\boldsymbol{\mathcal{Q}}\boldsymbol{\mathcal{B}}=\boldsymbol{0}$.
The idea is to separate the state $\boldsymbol{x}=\boldsymbol{\mathcal{P}}\boldsymbol{x}+\boldsymbol{\mathcal{Q}}\boldsymbol{x}$
and adjoint state $\boldsymbol{\lambda}=\boldsymbol{\mathcal{P}}^{T}\boldsymbol{\lambda}+\boldsymbol{\mathcal{Q}}^{T}\boldsymbol{\lambda}$
as well as their evolution equations with the help of $\boldsymbol{\mathcal{P}}$
and $\boldsymbol{\mathcal{Q}}$. The controlled state equation (\ref{eq:StateEq})
is split as
\begin{align}
\boldsymbol{\mathcal{P}}\boldsymbol{\dot{x}} & =\boldsymbol{\mathcal{P}}\boldsymbol{R}+\boldsymbol{\mathcal{B}}\boldsymbol{u}, & \boldsymbol{\mathcal{Q}}\boldsymbol{\dot{x}} & =\boldsymbol{\mathcal{Q}}\boldsymbol{R}.\label{eq:StateEqSeparated}
\end{align}
After multiplication with $\boldsymbol{\mathcal{B}}^{T}\boldsymbol{\mathcal{S}}$,
the first equation yields an expression for the control signal in
terms of $\boldsymbol{x}$, 
\begin{align}
\boldsymbol{u} & =\boldsymbol{\mathcal{B}}^{g}\left(\boldsymbol{\dot{x}}-\boldsymbol{R}\right), & \boldsymbol{\mathcal{B}}^{g} & =\left(\boldsymbol{\mathcal{B}}^{T}\boldsymbol{\mathcal{S}}\boldsymbol{\mathcal{B}}\right)^{-1}\boldsymbol{\mathcal{B}}^{T}\boldsymbol{\mathcal{S}}.\label{eq:ControlSignal}
\end{align}
Inserting $\boldsymbol{u}$ from Eq. (\ref{eq:ControlSignal}) in
Eq. (\ref{eq:StationarityCondition}), multiplying by $\boldsymbol{\mathcal{B}}^{gT}$
from the left and using $\boldsymbol{\mathcal{P}}=\boldsymbol{\mathcal{B}}\boldsymbol{\mathcal{B}}^{g}$
yields
\begin{align}
\boldsymbol{\mathcal{P}}^{T}\boldsymbol{\lambda} & =-\varepsilon^{2}\boldsymbol{\Gamma}\left(\boldsymbol{\dot{x}}-\boldsymbol{R}\right),\label{eq:Plambda}
\end{align}
with symmetric $n\times n$ matrix $\boldsymbol{\Gamma}=\boldsymbol{\mathcal{B}}^{gT}\boldsymbol{\mathcal{B}}^{g}$
of rank $p$. The adjoint state equation (\ref{eq:AdjointEq}) is
split analogously to Eq. (\ref{eq:StateEqSeparated}). Subsequently,
Eq. (\ref{eq:Plambda}) is used to eliminate $\boldsymbol{\mathcal{P}}^{T}\boldsymbol{\lambda}$
as well as $\boldsymbol{\mathcal{P}}^{T}\boldsymbol{\dot{\lambda}}$
from all equations. See the Supplemental Material (SM) \cite{supplement}
for a detailed derivation. The rearrangement results in the following,
singularly perturbed system of $2\left(n-p\right)$ linearly independent
ODEs and $p$ linearly independent second order differential equations,
\begin{align}
-\boldsymbol{\mathcal{Q}}^{T}\boldsymbol{\dot{\lambda}} & =\boldsymbol{\mathcal{Q}}^{T}\boldsymbol{w}_{\varepsilon}+\boldsymbol{\mathcal{Q}}^{T}\boldsymbol{\mathcal{S}}\boldsymbol{\mathcal{Q}}\left(\boldsymbol{x}-\boldsymbol{x}_{d}\right),\label{eq:Rearranged1}\\
\varepsilon^{2}\boldsymbol{\Gamma}\boldsymbol{\ddot{x}} & =\varepsilon^{2}\boldsymbol{\mathcal{P}}^{T}\left(\boldsymbol{\Gamma}\nabla\boldsymbol{R}\boldsymbol{\dot{x}}-\boldsymbol{\dot{\Gamma}}\left(\boldsymbol{\dot{x}}-\boldsymbol{R}\right)\right)+\boldsymbol{\mathcal{P}}^{T}\boldsymbol{w}_{\varepsilon}\nonumber \\
 & +\boldsymbol{\mathcal{P}}^{T}\boldsymbol{\mathcal{\dot{Q}}}^{T}\boldsymbol{\mathcal{Q}}^{T}\boldsymbol{\lambda}+\boldsymbol{\mathcal{P}}^{T}\boldsymbol{\mathcal{S}}\boldsymbol{\mathcal{P}}\left(\boldsymbol{x}-\boldsymbol{x}_{d}\right),\label{eq:Rearranged2}\\
\boldsymbol{\mathcal{Q}}\boldsymbol{\dot{x}} & =\boldsymbol{\mathcal{Q}}\boldsymbol{R}.\label{eq:Rearranged3}
\end{align}
We introduced the $n\times n$ matrix $\boldsymbol{\mathcal{W}}$
with entries $\mathcal{W}_{ij}=\sum_{k,l}\partial_{j}\mathcal{B}_{il}\mathcal{B}_{lk}^{g}\left(\dot{x}_{k}-R_{k}\right)$
and the vector
\begin{align}
\boldsymbol{w}_{\varepsilon} & =\left(\nabla\boldsymbol{R}^{T}+\boldsymbol{\mathcal{W}}^{T}\right)\left(\boldsymbol{\mathcal{Q}}^{T}\boldsymbol{\lambda}-\varepsilon^{2}\boldsymbol{\Gamma}\left(\boldsymbol{\dot{x}}-\boldsymbol{R}\right)\right).
\end{align}
We emphasize that the rearrangement requires no approximation and
is valid for all affine control systems with a cost functional of
the form Eq. (\ref{eq:OptimalTrajectoryTrackingFunctionalGeneral}).

\paragraph{Outer equations and linearizing assumption.}

The outer equations are obtained by setting $\varepsilon=0$ in Eqs.
(\ref{eq:Rearranged1})-(\ref{eq:Rearranged3}) and subsequently multiplying
Eq. (\ref{eq:Rearranged2}) with $\boldsymbol{\Omega}$ from the left.
Denoting the outer solutions by upper-case letters, $\boldsymbol{X}\left(t\right)=\boldsymbol{x}\left(t\right)$
and $\boldsymbol{\Lambda}\left(t\right)=\boldsymbol{\lambda}\left(t\right)$,
we obtain
\begin{align}
\boldsymbol{\mathcal{Q}}^{T}\boldsymbol{\dot{\Lambda}} & =-\boldsymbol{\mathcal{Q}}^{T}\boldsymbol{w}_{0}-\boldsymbol{\mathcal{Q}}^{T}\boldsymbol{\mathcal{S}}\boldsymbol{\mathcal{Q}}\left(\boldsymbol{X}-\boldsymbol{x}_{d}\right),\label{eq:OuterLeadingOrder1}\\
\boldsymbol{\mathcal{P}}\boldsymbol{X} & =\boldsymbol{\mathcal{P}}\boldsymbol{x}_{d}-\boldsymbol{\Omega}\left(\boldsymbol{w}_{0}+\boldsymbol{\mathcal{\dot{Q}}}^{T}\boldsymbol{\mathcal{Q}}^{T}\boldsymbol{\Lambda}\right),\label{eq:OuterLeadingOrder2}\\
\boldsymbol{\mathcal{Q}}\boldsymbol{\dot{X}} & =\boldsymbol{\mathcal{Q}}\boldsymbol{R},\label{eq:OuterLeadingOrder3}
\end{align}
with $\boldsymbol{w}_{0}=\left(\nabla\boldsymbol{R}^{T}+\boldsymbol{\mathcal{W}}^{T}\right)\boldsymbol{\mathcal{Q}}^{T}\boldsymbol{\Lambda}$
and $\boldsymbol{\mathcal{P}}^{T}\boldsymbol{\Lambda}=\boldsymbol{0}$.
In general, Eqs. (\ref{eq:OuterLeadingOrder1})-(\ref{eq:OuterLeadingOrder3})
are nonlinear because $\boldsymbol{\Omega}$, $\boldsymbol{\mathcal{P}}$,
$\boldsymbol{\mathcal{Q}}$ , $\boldsymbol{R}$, and $\boldsymbol{w}_{0}$
may all depend on $\boldsymbol{X}$.

The linearizing assumption consists of two parts. First, the matrix
$\boldsymbol{\Omega}\left(\boldsymbol{x}\right)$ is assumed to be
independent of $\boldsymbol{x}$. This assumption implies constant
projectors $\boldsymbol{\mathcal{P}}$ and $\boldsymbol{\mathcal{Q}}$
but not a constant input matrix $\boldsymbol{\mathcal{B}}\left(\boldsymbol{x}\right)$.
Second, the nonlinearity $\boldsymbol{R}\left(\boldsymbol{x}\right)$
is assumed to have the following structure with respect to the input
matrix, 
\begin{align}
\boldsymbol{\mathcal{Q}}\boldsymbol{R}\left(\boldsymbol{x}\right) & =\boldsymbol{\mathcal{Q}}\boldsymbol{\mathcal{A}}\boldsymbol{x}+\boldsymbol{\mathcal{Q}}\boldsymbol{b},\label{eq:LinearizingAssumptionS_2}
\end{align}
with constant $n\times n$ matrix $\boldsymbol{\mathcal{A}}$ and
constant $n$-component vector $\boldsymbol{b}$. Equation (\ref{eq:OuterLeadingOrder2})
together with the linearizing assumption yields an explicit expression
for $\boldsymbol{\mathcal{P}}\boldsymbol{X}$, 
\begin{align}
\boldsymbol{\mathcal{P}}\boldsymbol{X} & =\boldsymbol{\mathcal{P}}\boldsymbol{x}_{d}-\boldsymbol{\Omega}\boldsymbol{\mathcal{A}}^{T}\boldsymbol{\mathcal{Q}}^{T}\boldsymbol{\Lambda}.\label{eq:PxLinear}
\end{align}
Using Eqs. (\ref{eq:LinearizingAssumptionS_2}) and (\ref{eq:PxLinear})
in Eq. (\ref{eq:OuterLeadingOrder3}), we obtain $2\left(n-p\right)$
linearly independent ODEs for $\boldsymbol{\mathcal{Q}}^{T}\boldsymbol{\Lambda}$
and $\boldsymbol{\mathcal{Q}}\boldsymbol{X}$, 
\begin{align}
\negthickspace\left(\begin{array}{c}
\boldsymbol{\mathcal{Q}}^{T}\boldsymbol{\dot{\Lambda}}\\
\boldsymbol{\mathcal{Q}}\boldsymbol{\dot{X}}
\end{array}\right) & =\boldsymbol{\mathcal{M}}\left(\begin{array}{c}
\boldsymbol{\mathcal{Q}}^{T}\boldsymbol{\Lambda}\\
\boldsymbol{\mathcal{Q}}\boldsymbol{X}
\end{array}\right)+\left(\begin{array}{c}
\boldsymbol{\mathcal{Q}}^{T}\boldsymbol{\mathcal{S}}\boldsymbol{\mathcal{Q}}\boldsymbol{x}_{d}\\
\boldsymbol{\mathcal{Q}}\boldsymbol{\mathcal{A}}\boldsymbol{\mathcal{P}}\boldsymbol{x}_{d}+\boldsymbol{\mathcal{Q}}\boldsymbol{b}
\end{array}\right),\label{eq:OuterEquations}
\end{align}
with
\begin{align}
\boldsymbol{\mathcal{M}} & =\left(\begin{array}{cc}
-\boldsymbol{\mathcal{Q}}^{T}\boldsymbol{\mathcal{A}}^{T}\boldsymbol{\mathcal{Q}}^{T} & -\boldsymbol{\mathcal{Q}}^{T}\boldsymbol{\mathcal{S}}\boldsymbol{\mathcal{Q}}\\
-\boldsymbol{\mathcal{Q}}\boldsymbol{\mathcal{A}}\boldsymbol{\Omega}\boldsymbol{\mathcal{A}}^{T}\boldsymbol{\mathcal{Q}}^{T} & \boldsymbol{\mathcal{Q}}\boldsymbol{\mathcal{A}}\boldsymbol{\mathcal{Q}}
\end{array}\right).
\end{align}

\paragraph{Inner equations and matching.}

The initial inner equation, valid at the beginning, or left side,
of the time domain, is obtained by introducing functions $\boldsymbol{X}_{L}\left(\tau_{L}\right)=\boldsymbol{x}\left(t\right)$
and $\boldsymbol{\Lambda}_{L}\left(\tau_{L}\right)=\boldsymbol{\lambda}\left(t\right)$
with a rescaled time $\tau_{L}=\left(t-t_{0}\right)/\varepsilon\geq0$.
Performing a perturbation expansion of Eqs. (\ref{eq:Rearranged1})-(\ref{eq:Rearranged3})
up to leading order in $\varepsilon$ together with the linearizing
assumption yields an ODE for $\boldsymbol{\mathcal{P}}\boldsymbol{X}_{L}$,
\begin{align}
\left(\boldsymbol{\Gamma}\boldsymbol{X}_{L}'\right)' & =\boldsymbol{\mathcal{P}}^{T}\boldsymbol{\mathcal{A}}^{T}\boldsymbol{\mathcal{Q}}^{T}\boldsymbol{\Lambda}_{L}-\boldsymbol{\mathcal{P}}^{T}\boldsymbol{\mathcal{V}}^{T}\boldsymbol{\Gamma}\boldsymbol{X}_{L}'\nonumber \\
 & +\boldsymbol{\mathcal{P}}^{T}\boldsymbol{\mathcal{S}}\boldsymbol{\mathcal{P}}\left(\boldsymbol{X}_{L}-\boldsymbol{x}_{d}\left(t_{0}\right)\right),\label{eq:LeftInner}
\end{align}
and $\boldsymbol{\mathcal{P}}^{T}\boldsymbol{\Lambda}_{L}=\boldsymbol{0}$,
$\boldsymbol{\mathcal{Q}}^{T}\boldsymbol{\Lambda}_{L}'=\boldsymbol{0}$,
$\boldsymbol{\mathcal{Q}}\boldsymbol{X}_{L}'=\boldsymbol{0}$. We
introduced the $n\times n$ matrix $\boldsymbol{\mathcal{V}}$ with
entries $\mathcal{V}_{ij}=\sum_{k,l}\partial_{j}\mathcal{B}_{il}\mathcal{B}_{lk}^{g}X_{L,k}'$,
and $\left(\cdot\right)'=\frac{\partial}{\partial\tau_{L}}\left(\cdot\right)$.
An analogous procedure yields the right inner equations for $\boldsymbol{X}_{R}\left(\tau_{R}\right)$
and $\boldsymbol{\Lambda}_{R}\left(\tau_{R}\right)$ with rescaled
time $\tau_{R}=\left(t_{1}-t\right)/\varepsilon\geq0$ valid at the
right end of the time domain. They are identical in form to the left
inner equations. All inner and outer equations must be solved with
the matching conditions \cite{bender1999advanced},
\begin{flalign}
\lim_{\tau_{L}\rightarrow\infty}\boldsymbol{X}_{L}\left(\tau_{L}\right)=\boldsymbol{X}\left(t_{0}\right), & \lim_{\tau_{L}\rightarrow\infty}\boldsymbol{\Lambda}_{L}\left(\tau_{L}\right)=\boldsymbol{\Lambda}\left(t_{0}\right),\\
\lim_{\tau_{R}\rightarrow\infty}\boldsymbol{X}_{R}\left(\tau_{R}\right)=\boldsymbol{X}\left(t_{1}\right), & \lim_{\tau_{R}\rightarrow\infty}\boldsymbol{\Lambda}_{R}\left(\tau_{R}\right)=\boldsymbol{\Lambda}\left(t_{1}\right),
\end{flalign}
and the boundary conditions for the state, $\boldsymbol{X}_{L}\left(0\right)=\boldsymbol{x}_{0}$,
$\boldsymbol{X}_{R}\left(0\right)=\boldsymbol{x}_{1}$, see Eq. (\ref{eq:ControlledState}).
From the constancy of $\boldsymbol{\mathcal{Q}}\boldsymbol{X}_{L}$
and $\boldsymbol{\mathcal{Q}}\boldsymbol{X}_{R}$ follow immediately
the boundary conditions for the outer equations (\ref{eq:OuterEquations}),
$\boldsymbol{\mathcal{Q}}\boldsymbol{X}\left(t_{0}\right)=\boldsymbol{\mathcal{Q}}\boldsymbol{x}_{0}$
and $\boldsymbol{\mathcal{Q}}\boldsymbol{X}\left(t_{1}\right)=\boldsymbol{\mathcal{Q}}\boldsymbol{x}_{1}$.
The only non-constant inner solutions are $\boldsymbol{\mathcal{P}}\boldsymbol{X}_{L}$
and $\boldsymbol{\mathcal{P}}\boldsymbol{X}_{R}$. They provide a
connection from the the boundary conditions $\boldsymbol{\mathcal{P}}\boldsymbol{x}_{0}$
and $\boldsymbol{\mathcal{P}}\boldsymbol{x}_{1}$ to the outer solution
$\boldsymbol{\mathcal{P}}\boldsymbol{X}$ given by Eq. (\ref{eq:PxLinear}).
The connection is in form of a steep transition of width $\varepsilon$
known as a boundary layer \cite{bender1999advanced}. Note that the
inner equations may still be nonlinear due to a possible dependence
of $\boldsymbol{\mathcal{V}}$ and $\boldsymbol{\Gamma}$ on $\boldsymbol{X}_{L/R}$.
This nonlinearity originates from the input matrix $\boldsymbol{\mathcal{B}}$,
and vanishes for constant $\boldsymbol{\mathcal{B}}$. See the SM
\cite{supplement} for details of the derivation.\\
The approximate solution to the necessary optimality conditions depends
on the initial state $\boldsymbol{x}\left(t_{0}\right)=\boldsymbol{x}_{0}$.
The system can either be prepared in $\boldsymbol{x}_{0}$, or $\boldsymbol{x}_{0}$
is obtained by measurements. Once $\boldsymbol{x}_{0}$ is known,
no further information about the controlled system's state is necessary
to compute the control, resulting in an open-loop control. Measurements
$\boldsymbol{\tilde{x}}=\boldsymbol{x}\left(\tilde{t}\right)$ of
the state performed at a later time $\tilde{t}>t_{0}$ can be used
to update the control with $\boldsymbol{x}_{0}=\boldsymbol{\tilde{x}}$
as the new initial condition, leading to a sampled-data feedback law.
In the limit of continuous monitoring of the system's state $\boldsymbol{x}\left(t\right)$,
the initial conditions $\boldsymbol{x}_{0}=\boldsymbol{x}\left(t\right)$
become functions of the current state of the controlled system itself,
and the control becomes a continuous time feedback law \cite{bryson1969applied}\textit{.}
\begin{figure}
\includegraphics[width=8.6cm]{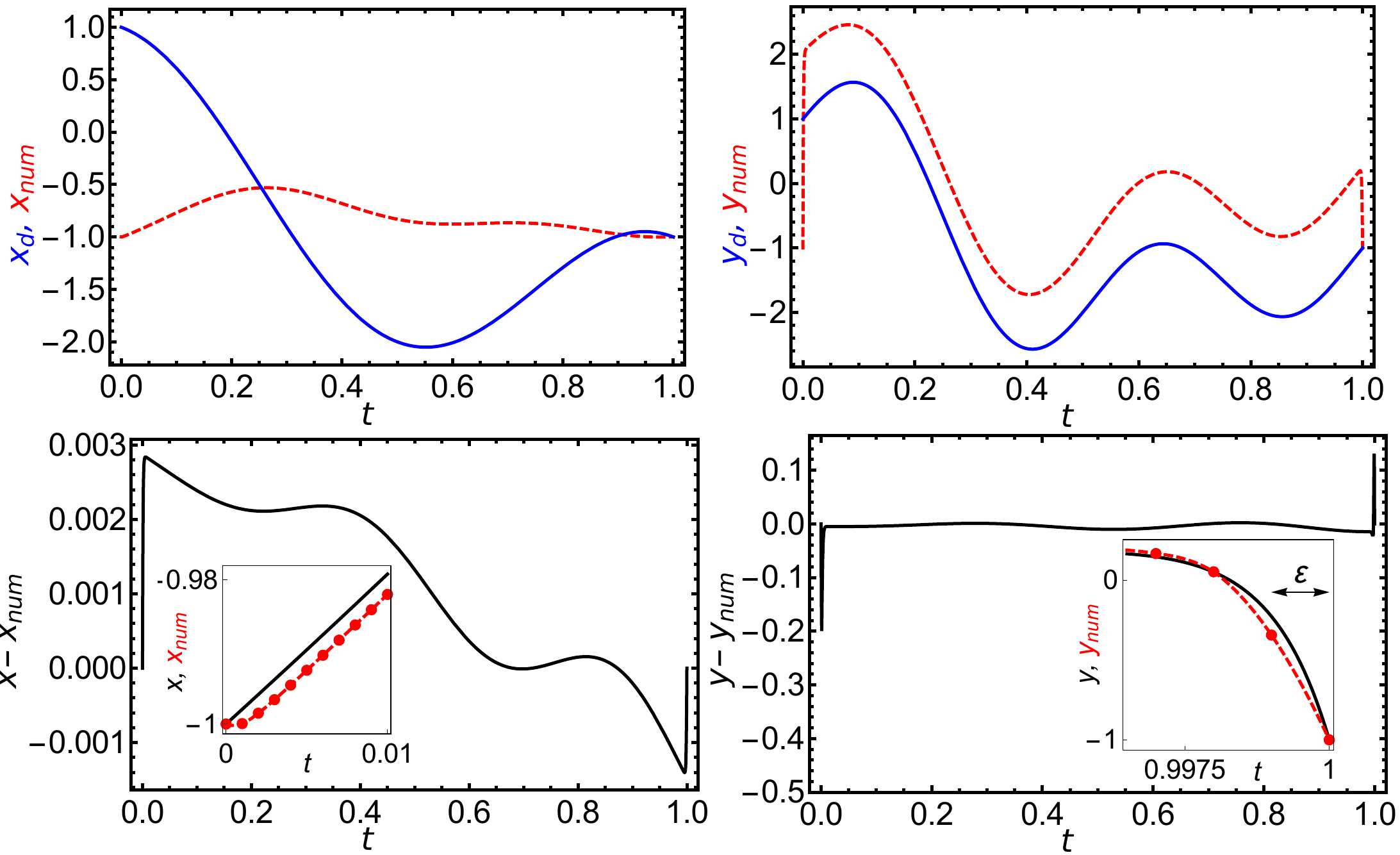}\protect\caption{\label{fig:Figure1}Optimal trajectory tracking for the damped mathematical
pendulum. Top: While the numerically obtained, optimally controlled
angular velocity $y_{num}$ (red dashed) resembles the desired velocity
$y_{d}$ (blue solid) in shape (right), the angular displacement $x_{num}$
(left) is far off its desired counterpart $x_{d}$. Bottom: The difference
between analytical and numerical solution for the optimally controlled
state reveals excellent agreement except close to the time domain
boundaries. To leading order in the small parameter $\varepsilon$,
the analytical results predicts a boundary layer of width $\varepsilon$
for $y$ (right inset) but not for $x$ (left inset). As indicated
by the red dots, the temporal resolution $\Delta t=\varepsilon=10^{-3}$
is just large enough to resolve the boundary layer.}
\end{figure}

\paragraph{Results for a two-dimensional dynamical system.}

We compare the approximate analytical solution for small but finite
$\varepsilon>0$ with numerical results for the optimally controlled
damped mathematical pendulum. The model is of the form Eq. (\ref{eq:xyState})
which satisfies the linearizing assumption. The EOM are 
\begin{align}
\dot{x} & =y, & \dot{y} & =-\frac{1}{2}y-\sin\left(x\right)+\left(1+\frac{1}{4}x^{2}\right)u,\label{eq:MathematicalPendulum}
\end{align}
with $x$ denoting the angular displacement and $y$ the angular velocity.
The desired trajectory $\boldsymbol{x}_{d}=\left(x_{d},y_{d}\right)^{T}$
is $x_{d}\left(t\right)=\cos\left(2\pi t\right)-2t$, $y_{d}\left(t\right)=x_{d}\left(t\right)+\sin\left(4\pi t\right)$.
The terminal conditions $x\left(1\right)=x_{d}\left(1\right)=-1$
and $y\left(1\right)=y_{d}\left(1\right)=-1$ are chosen to lie on
the desired trajectory, while the initial conditions $x\left(0\right)=-1$
and $y\left(0\right)=-1$ do not. Explicit analytical expressions
for the optimally controlled composite state solution and the optimal
control signal can be found in the SM \cite{supplement}. The prescribed
desired trajectory is not exactly realizable, which can be understood
from physical reasoning. A point mass governed by Newton's EOM can
only trace out desired trajectories in phase space, spanned by $x$
and $y$, which satisfy $\dot{x}_{d}=y_{d}$. This relation simply
defines the velocity, and no control force of whatever amplitude can
change that definition \cite{lober2016exactly}. Numerical computations
are performed with the ACADO Toolkit \cite{Houska2011a}, an open
source program package for solving optimal control problems. The problem
is solved on a time interval of length $1$ with a time step width
$\Delta t=10^{-3}$. A straightforward numerical solution of the necessary
optimality conditions Eqs. (\ref{eq:StationarityCondition})-(\ref{eq:AdjointEq})
is usually impossible due to the mixed boundary conditions. While
half of the $2n$ boundary conditions are imposed at the initial time,
the other half are imposed at the terminal time. This typically requires
iterative algorithms such as shooting, and renders optimal trajectory
tracking computationally expensive.\\
Figure \ref{fig:Figure1} top compares the desired trajectory $\boldsymbol{x}_{d}$
for angular displacement $x$ (left) and velocity $y$ (right) with
the numerically computed, optimally controlled state $\boldsymbol{x}$.
While the controlled velocity looks similar in shape to the desired
trajectory, the controlled displacement is way off. Even though the
terminal conditions comply with the desired trajectory, the solution
for $y$ exhibits some very steep transients at both ends of the time
interval. These transitions are the boundary layers of width $\varepsilon=10^{-3}$
described by the inner solutions. On this scale, there is no visible
difference between the analytically and numerically obtained controlled
state. We plot the difference between them in Fig. \ref{fig:Figure1}
bottom. Apart from the boundary layer region, the agreement is excellent.
Deviations in the boundary layer region can be understood from the
relatively poor numerical resolution equal to the boundary layer width,
$\Delta t=\varepsilon$, see the insets for closeups of the boundary
layers regions. Increasing $\Delta t$ leads to larger numerical errors
but decreasing $\Delta t$ quickly increases computation time. While
the analytical leading order result predicts a boundary layer for
$y$ but not for $x$, the numerical solution for $x$ exhibits a
tiny boundary layer, which we expect analytically to arise from higher
order contributions in the perturbation expansion.

\paragraph{Exact solution for $\varepsilon=0$.}

For $\varepsilon=0$, the inner solutions do not play a role because
the boundary layers exhibited by $\boldsymbol{\mathcal{P}}\boldsymbol{x}$
degenerate to jumps,
\begin{align}
\boldsymbol{\mathcal{P}}\boldsymbol{x} & =\begin{cases}
\boldsymbol{\mathcal{P}}\boldsymbol{x}_{0}, & t=t_{0},\\
\boldsymbol{\mathcal{P}}\boldsymbol{x}_{d}-\boldsymbol{\Omega}\boldsymbol{\mathcal{A}}^{T}\boldsymbol{\mathcal{Q}}^{T}\boldsymbol{\Lambda}, & t_{0}<t<t_{1},\\
\boldsymbol{\mathcal{P}}\boldsymbol{x}_{1}, & t=t_{1}.
\end{cases}\label{eq:PXSolution}
\end{align}
The jump heights are fully determined by the outer solution $\boldsymbol{\mathcal{Q}}^{T}\boldsymbol{\Lambda}$
and the boundary conditions. The parts $\boldsymbol{\mathcal{P}}^{T}\boldsymbol{\lambda}=\boldsymbol{\mathcal{P}}^{T}\boldsymbol{\Lambda}=\boldsymbol{0}$
vanish identically for all times. The parts $\boldsymbol{\mathcal{Q}}^{T}\boldsymbol{\lambda}=\boldsymbol{\mathcal{Q}}^{T}\boldsymbol{\Lambda}$
and $\boldsymbol{\mathcal{Q}}\boldsymbol{x}=\boldsymbol{\mathcal{Q}}\boldsymbol{X}$
are from the linear outer equations (\ref{eq:OuterEquations}). Cast
in terms of the state transition matrix $\boldsymbol{\Phi}\left(t,t_{0}\right)=\exp\left(\boldsymbol{\mathcal{M}}\left(t-t_{0}\right)\right)$,
the solutions become 
\begin{align}
\left(\begin{array}{c}
\boldsymbol{\mathcal{Q}}^{T}\boldsymbol{\Lambda}\left(t\right)\\
\boldsymbol{\mathcal{Q}}\boldsymbol{X}\left(t\right)
\end{array}\right) & =\boldsymbol{\Phi}\left(t,t_{0}\right)\left(\begin{array}{c}
\boldsymbol{\mathcal{Q}}^{T}\boldsymbol{\Lambda}_{\text{init}}\\
\boldsymbol{\mathcal{Q}}\boldsymbol{x}_{0}
\end{array}\right)\\
 & +\int_{t_{0}}^{t}d\tau\boldsymbol{\Phi}\left(t,\tau\right)\left(\begin{array}{c}
\boldsymbol{\mathcal{Q}}^{T}\boldsymbol{\mathcal{S}}\boldsymbol{\mathcal{Q}}\boldsymbol{x}_{d}\left(\tau\right)\\
\boldsymbol{\mathcal{Q}}\boldsymbol{\mathcal{A}}\boldsymbol{\mathcal{P}}\boldsymbol{x}_{d}\left(\tau\right)+\boldsymbol{\mathcal{Q}}\boldsymbol{b}
\end{array}\right).\nonumber 
\end{align}
The parameter $\boldsymbol{\mathcal{Q}}^{T}\boldsymbol{\Lambda}_{\text{init}}$
is determined from the terminal condition $\boldsymbol{\mathcal{Q}}\boldsymbol{x}\left(t_{1}\right)=\boldsymbol{\mathcal{Q}}\boldsymbol{x}_{1}$.
The control is given in terms of $\boldsymbol{x}$ by Eq. (\ref{eq:ControlSignal}).
We obtain (see the SM \cite{supplement} for a derivation),
\begin{align}
\boldsymbol{u} & =\begin{cases}
\boldsymbol{u}_{0}, & t=t_{0},\\
\boldsymbol{\mathcal{B}}^{g}\left(\boldsymbol{X}\right)\left(\boldsymbol{\dot{X}}-\boldsymbol{R}\left(\boldsymbol{X}\right)\right), & t_{0}<t<t_{1},\\
\boldsymbol{u}_{1}, & t=t_{1},
\end{cases}\label{eq:OptimalControl}
\end{align}
with
\begin{align}
\boldsymbol{u}_{i} & =\boldsymbol{\mathcal{B}}^{g}\left(\boldsymbol{x}_{i}\right)\left(\boldsymbol{\dot{X}}\left(t_{i}\right)-\boldsymbol{R}\left(\boldsymbol{x}_{i}\right)\right)\nonumber \\
 & +\left(-1\right)^{i}2\boldsymbol{\mathcal{B}}^{g}\left(\boldsymbol{x}_{i}\right)\left(\boldsymbol{X}\left(t_{i}\right)-\boldsymbol{x}_{i}\right)\delta\left(t-t_{i}\right).
\end{align}
Because of $\boldsymbol{u}\sim\boldsymbol{\mathcal{P}}\boldsymbol{\dot{x}}$,
the jumps of $\boldsymbol{\mathcal{P}}\boldsymbol{x}$ at the time
domain boundaries lead to divergences $\delta\left(t\right)$ in form
of Dirac delta functions located at the time domain boundaries.

\paragraph{Conclusions.}

For $\varepsilon>0$, the dynamics of an optimal control system takes
place in the in the combined state space of state $\boldsymbol{x}$
and co-state $\boldsymbol{\lambda}$ of dimension $2n$. For $\varepsilon=0$,
the dynamics is restricted by $2p$ algebraic equations, $\boldsymbol{\mathcal{P}}^{T}\boldsymbol{\lambda}=\boldsymbol{0}$
and Eq. (\ref{eq:PXSolution}), to a hypersurface of dimension $2\left(n-p\right)$,
also called a singular surface \cite{bryson1969applied}. At the initial
and terminal time, kicks in form of a Dirac delta function mediated
by control induce an instantaneous transition from $\boldsymbol{x}_{0}$
onto the singular surface and from the singular surface to $\boldsymbol{x}_{1}$,
respectively. These instantaneous transitions render the state components
$\boldsymbol{\mathcal{P}}\boldsymbol{x}$ discontinuous at the initial
and terminal time, respectively. For $\varepsilon>0$, the discontinuities
are smoothed out in form of boundary layers, i.e., continuous transition
regions with a slope and width controlled by $\varepsilon$. The control
signals are finite and exhibit a sharp peak at the time domain boundaries,
with an amplitude inversely proportional to $\varepsilon$. This general
picture of optimal trajectory tracking for small $\varepsilon$ is
independent of the linearizing assumption and remains true for arbitrary
affine control systems \cite{loeber2015optimal}.

In experiments and numerical simulations, it is impossible to generate
diverging control signals. In practice, a finite value $\varepsilon>0$
in Eq. (\ref{eq:OptimalTrajectoryTrackingFunctionalGeneral}) is indispensable.
Nevertheless, understanding the behavior of control systems in the
limit $\varepsilon\rightarrow0$ is very useful for applications.
For example, to prevent or at least minimize any steep transitions
and large control amplitudes, we may choose the initial state to lie
on the singular surface. Furthermore, a faithful numerical representation
of the solution to optimal control requires a temporal resolution
$\Delta t\lesssim\varepsilon$ to resolve the initial and terminal
boundary layers.

The linearizing assumption yields linear algebraic and differential
outer equations. While the inner equations may still be nonlinear
for $\varepsilon>0$, they reduce to jumps, with jump heights independent
of any details of the inner equations. Consequently, for $\varepsilon=0$,
the exact solution for $\boldsymbol{x}$, $\boldsymbol{\lambda}$,
and $\boldsymbol{u}$, is entirely given in terms of the linear outer
equations. For models of the form Eq. (\ref{eq:xyState}), this culminates
in the following, surprising conclusion. The controlled state becomes
independent of the external force $R$ as $\varepsilon=0$. Furthermore,
the dependence on the nonlinearity $b$ is restricted to the boundary
layers, which become irrelevant for $\varepsilon=0$. Consequently,
for a given desired trajectory $\boldsymbol{x}_{d}$, and given initial
and terminal conditions, and $\varepsilon=0$, all optimally controlled
states of Eq. (\ref{eq:xyState}) trace out the same trajectories
in phase space, independent of $R$ and $b$. However, note that the
control signal $u$ still depends on $R$ and $b$.

Intuitively, the result can be understood as follows. The state components
which and which are not acted upon directly by control are encoded
in the matrices $\boldsymbol{\mathcal{P}}$ and $\boldsymbol{\mathcal{Q}}$,
respectively. Assumption (\ref{eq:LinearizingAssumptionS_2}) of the
linearizing assumption demands that the control acts on the nonlinear
state equations, and all other equations are linear. For $\varepsilon=0$,
the control amplitude is unrestricted. The control dominates over
the nonlinearity and can absorb an arbitrarily large part $\boldsymbol{\mathcal{P}}\boldsymbol{R}$
of $\boldsymbol{R}$, resulting in linear evolution equations. The
linearizing assumption is satisfied by all models of the form $\dot{x}=a_{0}+a_{1}y+a_{2}x,\,\dot{y}=R\left(x,y\right)+b\left(x,y\right)u$,
which includes Eq. (\ref{eq:xyState}) and the FitzHugh-Nagumo model
with a control acting on the activator $y$ as special cases. Another
example is the SIR model for disease transmission with the transmission
rate serving as the control \cite{loeber2015optimal}.

Generalizations of the results presented here are possible. One example
are noisy controlled systems, for which fundamental results exist
for linear optimal control in form of the linear-quadratic regulator
\cite{bryson1969applied}. In this context, we mention Ref. \cite{kappen2005linear}
which presents a linear theory for the control of nonlinear stochastic
systems. However, the method in \cite{kappen2005linear} is restricted
to systems with identical numbers of control signals and state components,
$n=p$. Other possibilities are generalizations to spatio-temporal
systems as e.g. reaction-diffusion systems \cite{lober2014controlling}.
\begin{acknowledgments}
We thank the DFG for funding via GRK 1558.\vspace{-0.45cm}
\end{acknowledgments}

\bibliographystyle{apsrev4-1}
\bibliography{literature}
%%%%%%%%%% Merge with supplemental materials %%%%%%%%%% 
\onecolumngrid
\clearpage\begin{center}
\textbf{\large Supplemental Material: Linear structures in nonlinear optimal trajectory tracking}
\end{center}
\begin{center}
Jakob Löber*
\end{center}
\begin{center}
\textit{Institut für Theoretische Physik, EW 7-1, Hardenbergstraße 36, Technische Universität Berlin, 10623 Berlin}
\end{center}%%%%%%%%%% Prefix a "S" to all equations, figures, tables and reset the counter %%%%%%%%%%
\setcounter{equation}{0}
\setcounter{figure}{0}
\setcounter{table}{0}
\setcounter{page}{1}
%\makeatletter

\section*{\label{sec:AppendixA}Appendix A: Rearranging the necessary optimality
conditions}

\renewcommand{\theequation}{S\arabic{equation}}
\renewcommand{\theHequation}{S\arabic{equation}}
\renewcommand{\thefigure}{S\arabic{figure}}
\renewcommand{\theHfigure}{S\arabic{figure}}
\renewcommand{\bibnumfmt}[1]{[S#1]}
\renewcommand{\citenumfont}[1]{S#1}Here, we give a detailed version of the rearrangement of the necessary
optimality conditions,
\begin{align}
\boldsymbol{0} & =\varepsilon^{2}\boldsymbol{u}\left(t\right)+\boldsymbol{\mathcal{B}}^{T}\left(\boldsymbol{x}\left(t\right)\right)\boldsymbol{\lambda}\left(t\right),\label{eq:StationarityCondition-1}\\
\boldsymbol{\dot{x}}\left(t\right) & =\boldsymbol{R}\left(\boldsymbol{x}\left(t\right)\right)+\boldsymbol{\mathcal{B}}\left(\boldsymbol{x}\left(t\right)\right)\boldsymbol{u}\left(t\right),\label{eq:StateEq-1}\\
-\boldsymbol{\dot{\lambda}}\left(t\right) & =\left(\nabla\boldsymbol{R}^{T}\left(\boldsymbol{x}\left(t\right)\right)+\left(\nabla\boldsymbol{\mathcal{B}}\left(\boldsymbol{x}\left(t\right)\right)\boldsymbol{u}\left(t\right)\right)^{T}\right)\boldsymbol{\lambda}\left(t\right)+\boldsymbol{\mathcal{S}}\left(\boldsymbol{x}\left(t\right)-\boldsymbol{x}_{d}\left(t\right)\right),\label{eq:AdjointEq-1}
\end{align}
with the initial and terminal conditions
\begin{align}
\boldsymbol{x}\left(t_{0}\right) & =\boldsymbol{x}_{0}, & \boldsymbol{x}\left(t_{1}\right) & =\boldsymbol{x}_{1}.\label{eq:GenDynSysInitTermCond}
\end{align}
Together with the assumption $0\leq\varepsilon\ll1$, the rearrangement
will enable the reinterpretation of Eqs. (\ref{eq:StationarityCondition-1})-(\ref{eq:AdjointEq-1})
as a singularly perturbed system of differential equations. To shorten
the notation, the time argument of $\boldsymbol{x}=\boldsymbol{x}\left(t\right)$,
$\boldsymbol{\lambda}=\boldsymbol{\lambda}\left(t\right)$ and $\boldsymbol{u}=\boldsymbol{u}\left(t\right)$
is suppressed in the following and some abbreviating matrices are
introduced. Let the $n\times n$ matrix $\boldsymbol{\Omega}\left(\boldsymbol{x}\right)$
be defined by 
\begin{align}
\boldsymbol{\Omega}\left(\boldsymbol{x}\right) & =\boldsymbol{\mathcal{B}}\left(\boldsymbol{x}\right)\left(\boldsymbol{\mathcal{B}}^{T}\left(\boldsymbol{x}\right)\boldsymbol{\mathcal{S}}\boldsymbol{\mathcal{B}}\left(\boldsymbol{x}\right)\right)^{-1}\boldsymbol{\mathcal{B}}^{T}\left(\boldsymbol{x}\right).\label{eq:DefOmega}
\end{align}
A simple calculation shows that $\boldsymbol{\Omega}\left(\boldsymbol{x}\right)$
is symmetric, 
\begin{align}
\boldsymbol{\Omega}^{T}\left(\boldsymbol{x}\right) & =\boldsymbol{\mathcal{B}}\left(\boldsymbol{x}\right)\left(\boldsymbol{\mathcal{B}}^{T}\left(\boldsymbol{x}\right)\boldsymbol{\mathcal{S}}\boldsymbol{\mathcal{B}}\left(\boldsymbol{x}\right)\right)^{-T}\boldsymbol{\mathcal{B}}^{T}\left(\boldsymbol{x}\right)=\boldsymbol{\mathcal{B}}\left(\boldsymbol{x}\right)\left(\boldsymbol{\mathcal{B}}^{T}\left(\boldsymbol{x}\right)\boldsymbol{\mathcal{S}}\boldsymbol{\mathcal{B}}\left(\boldsymbol{x}\right)\right)^{-1}\boldsymbol{\mathcal{B}}^{T}\left(\boldsymbol{x}\right)=\boldsymbol{\Omega}\left(\boldsymbol{x}\right).
\end{align}
Note that $\boldsymbol{\mathcal{B}}^{T}\left(\boldsymbol{x}\right)\boldsymbol{\mathcal{S}}\boldsymbol{\mathcal{B}}\left(\boldsymbol{x}\right)=\left(\boldsymbol{\mathcal{B}}^{T}\left(\boldsymbol{x}\right)\boldsymbol{\mathcal{S}}\boldsymbol{\mathcal{B}}\left(\boldsymbol{x}\right)\right)^{T}$
is a symmetric $p\times p$ matrix because $\boldsymbol{\mathcal{S}}$
is symmetric by assumption, and the inverse of a symmetric matrix
is symmetric. Let the two $n\times n$ projection matrices $\boldsymbol{\mathcal{P}}\left(\boldsymbol{x}\right)$
and $\boldsymbol{\mathcal{Q}}\left(\boldsymbol{x}\right)$ be defined
by
\begin{align}
\boldsymbol{\mathcal{P}}\left(\boldsymbol{x}\right) & =\boldsymbol{\Omega}\left(\boldsymbol{x}\right)\boldsymbol{\mathcal{S}}, & \boldsymbol{\mathcal{Q}}\left(\boldsymbol{x}\right) & =\mathbf{1}-\boldsymbol{\mathcal{P}}\left(\boldsymbol{x}\right).
\end{align}
$\boldsymbol{\mathcal{P}}\left(\boldsymbol{x}\right)$ and $\boldsymbol{\mathcal{Q}}\left(\boldsymbol{x}\right)$
are idempotent, $\boldsymbol{\mathcal{P}}^{2}\left(\boldsymbol{x}\right)=\boldsymbol{\mathcal{P}}\left(\boldsymbol{x}\right)$
and $\boldsymbol{\mathcal{Q}}^{2}\left(\boldsymbol{x}\right)=\boldsymbol{\mathcal{Q}}\left(\boldsymbol{x}\right)$.
Furthermore, $\boldsymbol{\mathcal{P}}\left(\boldsymbol{x}\right)$
and $\boldsymbol{\mathcal{Q}}\left(\boldsymbol{x}\right)$ satisfy
the relations
\begin{align}
\boldsymbol{\mathcal{P}}\left(\boldsymbol{x}\right)\boldsymbol{\mathcal{B}}\left(\boldsymbol{x}\right) & =\boldsymbol{\mathcal{B}}\left(\boldsymbol{x}\right), & \boldsymbol{\mathcal{Q}}\left(\boldsymbol{x}\right)\boldsymbol{\mathcal{B}}\left(\boldsymbol{x}\right) & =\boldsymbol{0}, & \boldsymbol{\mathcal{B}}^{T}\left(\boldsymbol{x}\right)\boldsymbol{\mathcal{S}}\boldsymbol{\mathcal{P}}\left(\boldsymbol{x}\right) & =\boldsymbol{\mathcal{B}}^{T}\left(\boldsymbol{x}\right)\boldsymbol{\mathcal{S}}, & \boldsymbol{\mathcal{B}}^{T}\left(\boldsymbol{x}\right)\boldsymbol{\mathcal{S}}\boldsymbol{\mathcal{Q}} & =\boldsymbol{0}.
\end{align}
Computing the transposed of $\boldsymbol{\mathcal{P}}\left(\boldsymbol{x}\right)$
and $\boldsymbol{\mathcal{Q}}\left(\boldsymbol{x}\right)$ yields
\begin{align}
\boldsymbol{\mathcal{P}}^{T}\left(\boldsymbol{x}\right) & =\boldsymbol{\mathcal{S}}^{T}\boldsymbol{\Omega}^{T}\left(\boldsymbol{x}\right)=\boldsymbol{\mathcal{S}}\boldsymbol{\Omega}\left(\boldsymbol{x}\right)\neq\boldsymbol{\mathcal{P}}\left(\boldsymbol{x}\right),\label{eq:ProjectorSymmetry}
\end{align}
and analogously for $\boldsymbol{\mathcal{Q}}\left(\boldsymbol{x}\right)$.
Equation (\ref{eq:ProjectorSymmetry}) shows that $\boldsymbol{\mathcal{P}}\left(\boldsymbol{x}\right)$,
and therefore also $\boldsymbol{\mathcal{Q}}\left(\boldsymbol{x}\right)$,
is not symmetric. However, $\boldsymbol{\mathcal{P}}^{T}\left(\boldsymbol{x}\right)$
satisfies the convenient property 
\begin{align}
\boldsymbol{\mathcal{P}}^{T}\left(\boldsymbol{x}\right)\boldsymbol{\mathcal{S}} & =\boldsymbol{\mathcal{S}}\boldsymbol{\Omega}\left(\boldsymbol{x}\right)\boldsymbol{\mathcal{S}}=\boldsymbol{\mathcal{S}}\boldsymbol{\mathcal{P}}\left(\boldsymbol{x}\right),\label{eq:SPCommute}
\end{align}
which implies
\begin{align}
\boldsymbol{\mathcal{P}}^{T}\left(\boldsymbol{x}\right)\boldsymbol{\mathcal{S}} & =\boldsymbol{\mathcal{P}}^{T}\left(\boldsymbol{x}\right)\boldsymbol{\mathcal{P}}^{T}\left(\boldsymbol{x}\right)\boldsymbol{\mathcal{S}}=\boldsymbol{\mathcal{P}}^{T}\left(\boldsymbol{x}\right)\boldsymbol{\mathcal{S}}\boldsymbol{\mathcal{P}}\left(\boldsymbol{x}\right),\label{eq:SPCommute2}
\end{align}
and similarly for $\boldsymbol{\mathcal{S}}\boldsymbol{\mathcal{Q}}\left(\boldsymbol{x}\right)$.
The product of $\boldsymbol{\Omega}\left(\boldsymbol{x}\right)$ with
$\boldsymbol{\mathcal{P}}^{T}\left(\boldsymbol{x}\right)$ yields
\begin{align}
\boldsymbol{\Omega}\left(\boldsymbol{x}\right)\boldsymbol{\mathcal{P}}^{T}\left(\boldsymbol{x}\right) & =\boldsymbol{\Omega}\left(\boldsymbol{x}\right)\boldsymbol{\mathcal{S}}\boldsymbol{\Omega}\left(\boldsymbol{x}\right)=\boldsymbol{\mathcal{B}}\left(\boldsymbol{x}\right)\left(\boldsymbol{\mathcal{B}}^{T}\left(\boldsymbol{x}\right)\boldsymbol{\mathcal{S}}\boldsymbol{\mathcal{B}}\left(\boldsymbol{x}\right)\right)^{-1}\boldsymbol{\mathcal{B}}^{T}\left(\boldsymbol{x}\right)=\boldsymbol{\Omega}\left(\boldsymbol{x}\right),\label{eq:Eq4210}
\end{align}
and
\begin{align}
\boldsymbol{\Omega}\left(\boldsymbol{x}\right)\boldsymbol{\mathcal{P}}^{T}\left(\boldsymbol{x}\right)\boldsymbol{\mathcal{S}} & =\boldsymbol{\Omega}\left(\boldsymbol{x}\right)\boldsymbol{\mathcal{S}}=\boldsymbol{\mathcal{P}}\left(\boldsymbol{x}\right).\label{eq:Eq4211}
\end{align}
Let the $n\times n$ matrix $\boldsymbol{\Gamma}\left(\boldsymbol{x}\right)$
be defined by 
\begin{align}
\boldsymbol{\Gamma}\left(\boldsymbol{x}\right) & =\boldsymbol{\mathcal{S}}\boldsymbol{\mathcal{B}}\left(\boldsymbol{x}\right)\left(\boldsymbol{\mathcal{B}}^{T}\left(\boldsymbol{x}\right)\boldsymbol{\mathcal{S}}\boldsymbol{\mathcal{B}}\left(\boldsymbol{x}\right)\right)^{-2}\boldsymbol{\mathcal{B}}^{T}\left(\boldsymbol{x}\right)\boldsymbol{\mathcal{S}}.\label{eq:DefGamma}
\end{align}
$\boldsymbol{\Gamma}\left(\boldsymbol{x}\right)$ is symmetric, 
\begin{align}
\boldsymbol{\Gamma}^{T}\left(\boldsymbol{x}\right) & =\left(\boldsymbol{\mathcal{S}}\boldsymbol{\mathcal{B}}\left(\boldsymbol{x}\right)\left(\boldsymbol{\mathcal{B}}^{T}\left(\boldsymbol{x}\right)\boldsymbol{\mathcal{S}}\boldsymbol{\mathcal{B}}\left(\boldsymbol{x}\right)\right)^{-2}\boldsymbol{\mathcal{B}}^{T}\left(\boldsymbol{x}\right)\boldsymbol{\mathcal{S}}\right)^{T}=\boldsymbol{\mathcal{S}}\boldsymbol{\mathcal{B}}\left(\boldsymbol{x}\right)\left(\boldsymbol{\mathcal{B}}^{T}\left(\boldsymbol{x}\right)\boldsymbol{\mathcal{S}}\boldsymbol{\mathcal{B}}\left(\boldsymbol{x}\right)\right)^{-2T}\boldsymbol{\mathcal{B}}^{T}\left(\boldsymbol{x}\right)\boldsymbol{\mathcal{S}}=\boldsymbol{\Gamma}\left(\boldsymbol{x}\right),\label{eq:GammaSSymmetry}
\end{align}
and satisfies
\begin{align}
\boldsymbol{\Gamma}\left(\boldsymbol{x}\right)\boldsymbol{\mathcal{P}}\left(\boldsymbol{x}\right)= & \boldsymbol{\mathcal{S}}\boldsymbol{\mathcal{B}}\left(\boldsymbol{x}\right)\left(\boldsymbol{\mathcal{B}}^{T}\left(\boldsymbol{x}\right)\boldsymbol{\mathcal{S}}\boldsymbol{\mathcal{B}}\left(\boldsymbol{x}\right)\right)^{-2}\boldsymbol{\mathcal{B}}^{T}\left(\boldsymbol{x}\right)\boldsymbol{\mathcal{S}}\boldsymbol{\mathcal{B}}\left(\boldsymbol{x}\right)\left(\boldsymbol{\mathcal{B}}^{T}\left(\boldsymbol{x}\right)\boldsymbol{\mathcal{S}}\boldsymbol{\mathcal{B}}\left(\boldsymbol{x}\right)\right)^{-1}\boldsymbol{\mathcal{B}}^{T}\left(\boldsymbol{x}\right)\boldsymbol{\mathcal{S}}\nonumber \\
= & \boldsymbol{\mathcal{S}}\boldsymbol{\mathcal{B}}\left(\boldsymbol{x}\right)\left(\boldsymbol{\mathcal{B}}^{T}\left(\boldsymbol{x}\right)\boldsymbol{\mathcal{S}}\boldsymbol{\mathcal{B}}\left(\boldsymbol{x}\right)\right)^{-2}\boldsymbol{\mathcal{B}}^{T}\left(\boldsymbol{x}\right)\boldsymbol{\mathcal{S}}=\boldsymbol{\Gamma}\left(\boldsymbol{x}\right).
\end{align}
Transposing yields
\begin{align}
\left(\boldsymbol{\Gamma}\left(\boldsymbol{x}\right)\boldsymbol{\mathcal{P}}\left(\boldsymbol{x}\right)\right)^{T} & =\boldsymbol{\mathcal{P}}^{T}\left(\boldsymbol{x}\right)\boldsymbol{\Gamma}^{T}\left(\boldsymbol{x}\right)=\boldsymbol{\mathcal{P}}^{T}\left(\boldsymbol{x}\right)\boldsymbol{\Gamma}\left(\boldsymbol{x}\right)=\boldsymbol{\Gamma}\left(\boldsymbol{x}\right).\label{eq:PSGammaS}
\end{align}
The projectors $\boldsymbol{\mathcal{P}}\left(\boldsymbol{x}\right)$
and $\boldsymbol{\mathcal{Q}}\left(\boldsymbol{x}\right)$ are used
to partition the state $\boldsymbol{x}$ as 
\begin{align}
\boldsymbol{x} & =\boldsymbol{\mathcal{P}}\left(\boldsymbol{x}\right)\boldsymbol{x}+\boldsymbol{\mathcal{Q}}\left(\boldsymbol{x}\right)\boldsymbol{x}.\label{eq:PartitionState}
\end{align}
The controlled state equation (\ref{eq:StateEq-1}) is split in two
parts by multiplying with $\boldsymbol{\mathcal{P}}\left(\boldsymbol{x}\right)$
and with $\boldsymbol{\mathcal{Q}}\left(\boldsymbol{x}\right)$ from
the left,
\begin{align}
\boldsymbol{\mathcal{P}}\left(\boldsymbol{x}\right)\boldsymbol{\dot{x}} & =\boldsymbol{\mathcal{P}}\left(\boldsymbol{x}\right)\boldsymbol{R}\left(\boldsymbol{x}\right)+\boldsymbol{\mathcal{B}}\left(\boldsymbol{x}\right)\boldsymbol{u}, & \boldsymbol{\mathcal{Q}}\left(\boldsymbol{x}\right)\boldsymbol{\dot{x}} & =\boldsymbol{\mathcal{Q}}\left(\boldsymbol{x}\right)\boldsymbol{R}\left(\boldsymbol{x}\right).\label{eq:StatePartitionPQ}
\end{align}
The initial and terminal conditions are split up as well,
\begin{align}
\boldsymbol{\mathcal{P}}\left(\boldsymbol{x}\left(t_{0}\right)\right)\boldsymbol{x}\left(t_{0}\right) & =\boldsymbol{\mathcal{P}}\left(\boldsymbol{x}_{0}\right)\boldsymbol{x}_{0}, & \boldsymbol{\mathcal{Q}}\left(\boldsymbol{x}\left(t_{0}\right)\right)\boldsymbol{x}\left(t_{0}\right) & =\boldsymbol{\mathcal{Q}}\left(\boldsymbol{x}_{0}\right)\boldsymbol{x}_{0},\label{eq:PSQSInitCond}\\
\boldsymbol{\mathcal{P}}\left(\boldsymbol{x}\left(t_{1}\right)\right)\boldsymbol{x}\left(t_{1}\right) & =\boldsymbol{\mathcal{P}}\left(\boldsymbol{x}_{1}\right)\boldsymbol{x}_{1}, & \boldsymbol{\mathcal{Q}}\left(\boldsymbol{x}\left(t_{1}\right)\right)\boldsymbol{x}\left(t_{1}\right) & =\boldsymbol{\mathcal{Q}}\left(\boldsymbol{x}_{1}\right)\boldsymbol{x}_{1}.\label{eq:PSQSTermCond}
\end{align}
Multiplying the first equation of Eq. (\ref{eq:StatePartitionPQ})
by $\boldsymbol{\mathcal{B}}^{T}\left(\boldsymbol{x}\right)\boldsymbol{\mathcal{S}}$
from the left and using $\boldsymbol{\mathcal{B}}^{T}\left(\boldsymbol{x}\right)\boldsymbol{\mathcal{S}}\boldsymbol{\mathcal{P}}\left(\boldsymbol{x}\right)=\boldsymbol{\mathcal{B}}^{T}\left(\boldsymbol{x}\right)\boldsymbol{\mathcal{S}}$
yields an expression for the control $\boldsymbol{u}$ in terms of
the controlled state trajectory $\boldsymbol{x}$, 
\begin{align}
\boldsymbol{u} & =\left(\boldsymbol{\mathcal{B}}^{T}\left(\boldsymbol{x}\right)\boldsymbol{\mathcal{S}}\boldsymbol{\mathcal{B}}\left(\boldsymbol{x}\right)\right)^{-1}\boldsymbol{\mathcal{B}}^{T}\left(\boldsymbol{x}\right)\boldsymbol{\mathcal{S}}\left(\boldsymbol{\dot{x}}-\boldsymbol{R}\left(\boldsymbol{x}\right)\right)=\boldsymbol{\mathcal{B}}^{g}\left(\boldsymbol{x}\right)\left(\boldsymbol{\dot{x}}-\boldsymbol{R}\left(\boldsymbol{x}\right)\right).
\end{align}
The $p\times n$ matrix $\boldsymbol{\mathcal{B}}^{g}\left(\boldsymbol{x}\right)$
is defined by
\begin{align}
\boldsymbol{\mathcal{B}}^{g}\left(\boldsymbol{x}\right) & =\left(\boldsymbol{\mathcal{B}}^{T}\left(\boldsymbol{x}\right)\boldsymbol{\mathcal{S}}\boldsymbol{\mathcal{B}}\left(\boldsymbol{x}\right)\right)^{-1}\boldsymbol{\mathcal{B}}^{T}\left(\boldsymbol{x}\right)\boldsymbol{\mathcal{S}}.
\end{align}
The matrix $\boldsymbol{\mathcal{B}}^{g}\left(\boldsymbol{x}\right)$
can be used to rewrite the matrices $\boldsymbol{\mathcal{P}}\left(\boldsymbol{x}\right)$
and $\boldsymbol{\Gamma}\left(\boldsymbol{x}\right)$ as
\begin{align}
\boldsymbol{\mathcal{P}}\left(\boldsymbol{x}\right) & =\boldsymbol{\mathcal{B}}\left(\boldsymbol{x}\right)\boldsymbol{\mathcal{B}}^{g}\left(\boldsymbol{x}\right), & \boldsymbol{\Gamma}\left(\boldsymbol{x}\right) & =\boldsymbol{\mathcal{B}}^{gT}\left(\boldsymbol{x}\right)\boldsymbol{\mathcal{B}}^{g}\left(\boldsymbol{x}\right),\label{eq:Eq4222}
\end{align}
respectively.\\
The solution for $\boldsymbol{u}$ is inserted in the stationarity
condition Eq. (\ref{eq:StationarityCondition-1}) to yield
\begin{align}
\boldsymbol{0} & =\varepsilon^{2}\boldsymbol{u}^{T}+\boldsymbol{\lambda}^{T}\boldsymbol{\mathcal{B}}\left(\boldsymbol{x}\right)=\varepsilon^{2}\left(\boldsymbol{\dot{x}}^{T}-\boldsymbol{R}^{T}\left(\boldsymbol{x}\right)\right)\boldsymbol{\mathcal{B}}^{gT}\left(\boldsymbol{x}\right)+\boldsymbol{\lambda}^{T}\boldsymbol{\mathcal{B}}\left(\boldsymbol{x}\right).\label{eq:StationarityCondition_1}
\end{align}
Equation (\ref{eq:StationarityCondition_1}) is utilized to eliminate
any occurrence of the part $\boldsymbol{\mathcal{P}}^{T}\left(\boldsymbol{x}\right)\boldsymbol{\lambda}$
in all equations. In contrast to the state $\boldsymbol{x}$, cf.
Eq. (\ref{eq:PartitionState}), the co-state is split up with the
transposed projectors $\boldsymbol{\mathcal{P}}^{T}\left(\boldsymbol{x}\right)$
and $\boldsymbol{\mathcal{Q}}^{T}\left(\boldsymbol{x}\right)$,
\begin{align}
\boldsymbol{\lambda} & =\boldsymbol{\mathcal{P}}^{T}\left(\boldsymbol{x}\right)\boldsymbol{\lambda}+\boldsymbol{\mathcal{Q}}^{T}\left(\boldsymbol{x}\right)\boldsymbol{\lambda}.\label{eq:PartitionCoState}
\end{align}
Multiplying Eq. (\ref{eq:StationarityCondition_1}) with $\boldsymbol{\mathcal{B}}^{g}\left(\boldsymbol{x}\right)$
from the right and using Eq. (\ref{eq:Eq4222}) yields an expression
for $\boldsymbol{\mathcal{P}}^{T}\left(\boldsymbol{x}\right)\boldsymbol{\lambda}$,
\begin{align}
\boldsymbol{0} & =\varepsilon^{2}\left(\boldsymbol{\dot{x}}^{T}-\boldsymbol{R}^{T}\left(\boldsymbol{x}\right)\right)\boldsymbol{\Gamma}\left(\boldsymbol{x}\right)+\boldsymbol{\lambda}^{T}\boldsymbol{\mathcal{P}}\left(\boldsymbol{x}\right).
\end{align}
Transposing the last equation and exploiting the symmetry of $\boldsymbol{\Gamma}\left(\boldsymbol{x}\right)$,
Eq. (\ref{eq:GammaSSymmetry}), yields 
\begin{align}
\boldsymbol{\mathcal{P}}^{T}\left(\boldsymbol{x}\right)\boldsymbol{\lambda} & =-\varepsilon^{2}\boldsymbol{\Gamma}\left(\boldsymbol{x}\right)\left(\boldsymbol{\dot{x}}-\boldsymbol{R}\left(\boldsymbol{x}\right)\right).\label{eq:PSLambda}
\end{align}
Equation (\ref{eq:PSLambda}) is valid for all times $t_{0}\leq t\leq t_{1}$
such that we can apply the time derivative to get
\begin{align}
\boldsymbol{0} & =\varepsilon^{2}\boldsymbol{\Gamma}\left(\boldsymbol{x}\right)\left(\boldsymbol{\ddot{x}}-\nabla\boldsymbol{R}\left(\boldsymbol{x}\right)\boldsymbol{\dot{x}}\right)+\varepsilon^{2}\boldsymbol{\dot{\Gamma}}\left(\boldsymbol{x}\right)\left(\boldsymbol{\dot{x}}-\boldsymbol{R}\left(\boldsymbol{x}\right)\right)+\boldsymbol{\mathcal{\dot{P}}}^{T}\left(\boldsymbol{x}\right)\boldsymbol{\lambda}+\boldsymbol{\mathcal{P}}^{T}\left(\boldsymbol{x}\right)\boldsymbol{\dot{\lambda}}.\label{eq:PSLambdaDot_11}
\end{align}
The short hand notations
\begin{align}
\boldsymbol{\dot{\Gamma}}\left(\boldsymbol{x}\right) & =\dfrac{d}{dt}\boldsymbol{\Gamma}\left(\boldsymbol{x}\right)=\sum_{j=1}^{n}\dfrac{\partial}{\partial x_{j}}\boldsymbol{\Gamma}\left(\boldsymbol{x}\right)\dot{x}_{j}, & \boldsymbol{\mathcal{\dot{P}}}^{T}\left(\boldsymbol{x}\right) & =\dfrac{d}{dt}\boldsymbol{\mathcal{P}}^{T}\left(\boldsymbol{x}\right)=\sum_{j=1}^{n}\dfrac{\partial}{\partial x_{j}}\boldsymbol{\mathcal{P}}^{T}\left(\boldsymbol{x}\right)\dot{x}_{j},
\end{align}
were introduced in Eq. (\ref{eq:PSLambdaDot_11}). Splitting the co-state
$\boldsymbol{\lambda}$ as in Eq. (\ref{eq:PartitionCoState}) and
using Eq. (\ref{eq:PSLambda}) to eliminate $\boldsymbol{\mathcal{P}}^{T}\left(\boldsymbol{x}\right)\boldsymbol{\lambda}$
leads to 
\begin{align}
-\boldsymbol{\mathcal{P}}^{T}\left(\boldsymbol{x}\right)\boldsymbol{\dot{\lambda}} & =\varepsilon^{2}\boldsymbol{\Gamma}\left(\boldsymbol{x}\right)\left(\boldsymbol{\ddot{x}}-\nabla\boldsymbol{R}\left(\boldsymbol{x}\right)\boldsymbol{\dot{x}}\right)+\boldsymbol{\mathcal{\dot{P}}}^{T}\left(\boldsymbol{x}\right)\boldsymbol{\mathcal{Q}}^{T}\left(\boldsymbol{x}\right)\boldsymbol{\lambda}+\varepsilon^{2}\left(\boldsymbol{\dot{\Gamma}}\left(\boldsymbol{x}\right)-\boldsymbol{\mathcal{\dot{P}}}^{T}\left(\boldsymbol{x}\right)\boldsymbol{\Gamma}\left(\boldsymbol{x}\right)\right)\left(\boldsymbol{\dot{x}}-\boldsymbol{R}\left(\boldsymbol{x}\right)\right).\label{eq:PSLambdaDot_12}
\end{align}
Equation (\ref{eq:PSLambdaDot_12}) is an expression for $\boldsymbol{\mathcal{P}}^{T}\left(\boldsymbol{x}\right)\boldsymbol{\dot{\lambda}}$
independent of $\boldsymbol{\mathcal{P}}^{T}\left(\boldsymbol{x}\right)\boldsymbol{\lambda}$.\\
A similar procedure is performed for the adjoint equation (\ref{eq:AdjointEq-1}).
Eliminating the control signal $\boldsymbol{u}$ from Eq. (\ref{eq:AdjointEq-1})
gives
\begin{align}
-\boldsymbol{\dot{\lambda}} & =\left(\nabla\boldsymbol{R}^{T}\left(\boldsymbol{x}\right)+\boldsymbol{\mathcal{W}}^{T}\left(\boldsymbol{x},\boldsymbol{\dot{x}}\right)\right)\boldsymbol{\lambda}+\boldsymbol{\mathcal{S}}\left(\boldsymbol{x}-\boldsymbol{x}_{d}\left(t\right)\right).\label{eq:Eq4230-1}
\end{align}
Tho shorten the notation, the $n\times n$ matrix
\begin{align}
\boldsymbol{\mathcal{W}}\left(\boldsymbol{x},\boldsymbol{y}\right) & =\nabla\boldsymbol{\mathcal{B}}\left(\boldsymbol{x}\right)\boldsymbol{\mathcal{B}}^{g}\left(\boldsymbol{x}\right)\left(\boldsymbol{y}-\boldsymbol{R}\left(\boldsymbol{x}\right)\right)
\end{align}
with entries
\begin{align}
\mathcal{W}_{ij}\left(\boldsymbol{x},\boldsymbol{y}\right) & =\sum_{k=1}^{n}\sum_{l=1}^{p}\dfrac{\partial}{\partial x_{j}}\mathcal{B}_{il}\left(\boldsymbol{x}\right)\mathcal{B}_{lk}^{g}\left(\boldsymbol{x}\right)\left(y_{k}-R_{k}\left(\boldsymbol{x}\right)\right)\label{eq:DefW}
\end{align}
was introduced. With the help of the projectors $\boldsymbol{\mathcal{P}}^{T}$
and $\boldsymbol{\mathcal{Q}}^{T}$, Eq. (\ref{eq:Eq4230-1}) is split
up in two parts,
\begin{align}
-\boldsymbol{\mathcal{P}}^{T}\left(\boldsymbol{x}\right)\boldsymbol{\dot{\lambda}} & =\boldsymbol{\mathcal{P}}^{T}\left(\boldsymbol{x}\right)\left(\nabla\boldsymbol{R}^{T}\left(\boldsymbol{x}\right)+\boldsymbol{\mathcal{W}}^{T}\left(\boldsymbol{x},\boldsymbol{\dot{x}}\right)\right)\boldsymbol{\lambda}+\boldsymbol{\mathcal{P}}^{T}\left(\boldsymbol{x}\right)\boldsymbol{\mathcal{S}}\left(\boldsymbol{x}-\boldsymbol{x}_{d}\left(t\right)\right),\label{eq:PSLambdaDot_21}\\
-\boldsymbol{\mathcal{Q}}^{T}\left(\boldsymbol{x}\right)\boldsymbol{\dot{\lambda}} & =\boldsymbol{\mathcal{Q}}^{T}\left(\boldsymbol{x}\right)\left(\nabla\boldsymbol{R}^{T}\left(\boldsymbol{x}\right)+\boldsymbol{\mathcal{W}}^{T}\left(\boldsymbol{x},\boldsymbol{\dot{x}}\right)\right)\boldsymbol{\lambda}+\boldsymbol{\mathcal{Q}}^{T}\left(\boldsymbol{x}\right)\boldsymbol{\mathcal{S}}\left(\boldsymbol{x}-\boldsymbol{x}_{d}\left(t\right)\right).\label{eq:QSLambdaDot_21}
\end{align}
Using Eq. (\ref{eq:PSLambda}) to eliminate $\boldsymbol{\mathcal{P}}^{T}\left(\boldsymbol{x}\right)\boldsymbol{\lambda}$
in Eqs. (\ref{eq:PSLambdaDot_21}) and (\ref{eq:QSLambdaDot_21})
results in 
\begin{align}
-\boldsymbol{\mathcal{P}}^{T}\left(\boldsymbol{x}\right)\boldsymbol{\dot{\lambda}} & =\boldsymbol{\mathcal{P}}^{T}\left(\boldsymbol{x}\right)\boldsymbol{w}_{\varepsilon}\left(\boldsymbol{x},\boldsymbol{\dot{x}},\boldsymbol{\lambda}\right)+\boldsymbol{\mathcal{P}}^{T}\left(\boldsymbol{x}\right)\boldsymbol{\mathcal{S}}\left(\boldsymbol{x}-\boldsymbol{x}_{d}\left(t\right)\right),\label{eq:PSLambdaDot_22}\\
-\boldsymbol{\mathcal{Q}}^{T}\left(\boldsymbol{x}\right)\boldsymbol{\dot{\lambda}} & =\boldsymbol{\mathcal{Q}}^{T}\left(\boldsymbol{x}\right)\boldsymbol{w}_{\varepsilon}\left(\boldsymbol{x},\boldsymbol{\dot{x}},\boldsymbol{\lambda}\right)+\boldsymbol{\mathcal{Q}}^{T}\left(\boldsymbol{x}\right)\boldsymbol{\mathcal{S}}\left(\boldsymbol{x}-\boldsymbol{x}_{d}\left(t\right)\right).\label{eq:eq:QSLambdaDot_22}
\end{align}
with the abbreviation $\boldsymbol{w}_{\varepsilon}\left(\boldsymbol{x},\boldsymbol{y},\boldsymbol{z}\right)$
defined as the $n\times1$ vector
\begin{align}
\boldsymbol{w}_{\varepsilon}\left(\boldsymbol{x},\boldsymbol{y},\boldsymbol{z}\right) & =\left(\nabla\boldsymbol{R}^{T}\left(\boldsymbol{x}\right)+\boldsymbol{\mathcal{W}}^{T}\left(\boldsymbol{x},\boldsymbol{y}\right)\right)\left(\boldsymbol{\mathcal{Q}}^{T}\left(\boldsymbol{x}\right)\boldsymbol{z}-\varepsilon^{2}\boldsymbol{\Gamma}\left(\boldsymbol{x}\right)\left(\boldsymbol{y}-\boldsymbol{R}\left(\boldsymbol{x}\right)\right)\right).\label{eq:Defw}
\end{align}
Equations (\ref{eq:PSLambdaDot_12}) and (\ref{eq:PSLambdaDot_22})
are two independent expressions for $\boldsymbol{\mathcal{P}}^{T}\left(\boldsymbol{x}\right)\boldsymbol{\dot{\lambda}}$.
Combining them yields a second order differential equation independent
of $\boldsymbol{\mathcal{P}}^{T}\left(\boldsymbol{x}\right)\boldsymbol{\lambda}$
and $\boldsymbol{\mathcal{P}}^{T}\left(\boldsymbol{x}\right)\boldsymbol{\dot{\lambda}}$,
\begin{align}
\varepsilon^{2}\boldsymbol{\Gamma}\left(\boldsymbol{x}\right)\boldsymbol{\ddot{x}} & =\varepsilon^{2}\boldsymbol{\Gamma}\left(\boldsymbol{x}\right)\nabla\boldsymbol{R}\left(\boldsymbol{x}\right)\boldsymbol{\dot{x}}-\varepsilon^{2}\left(\boldsymbol{\dot{\Gamma}}\left(\boldsymbol{x}\right)-\boldsymbol{\mathcal{\dot{P}}}^{T}\left(\boldsymbol{x}\right)\boldsymbol{\Gamma}\left(\boldsymbol{x}\right)\right)\left(\boldsymbol{\dot{x}}-\boldsymbol{R}\left(\boldsymbol{x}\right)\right)\nonumber \\
 & +\boldsymbol{\mathcal{P}}^{T}\left(\boldsymbol{x}\right)\boldsymbol{w}_{\varepsilon}\left(\boldsymbol{x},\boldsymbol{\dot{x}},\boldsymbol{\lambda}\right)-\boldsymbol{\mathcal{\dot{P}}}^{T}\left(\boldsymbol{x}\right)\boldsymbol{\mathcal{Q}}^{T}\left(\boldsymbol{x}\right)\boldsymbol{\lambda}+\boldsymbol{\mathcal{P}}^{T}\left(\boldsymbol{x}\right)\boldsymbol{\mathcal{S}}\left(\boldsymbol{x}-\boldsymbol{x}_{d}\left(t\right)\right).\label{eq:PSLambdaDot_3}
\end{align}
Equation (\ref{eq:PSLambdaDot_3}) contains several time dependent
matrices which can be simplified. From Eq. (\ref{eq:PSGammaS}) follows
for the time derivative of $\boldsymbol{\Gamma}\left(\boldsymbol{x}\right)$
\begin{align}
\boldsymbol{\dot{\Gamma}}\left(\boldsymbol{x}\right) & =\boldsymbol{\mathcal{\dot{P}}}^{T}\left(\boldsymbol{x}\right)\boldsymbol{\Gamma}\left(\boldsymbol{x}\right)+\boldsymbol{\mathcal{P}}^{T}\left(\boldsymbol{x}\right)\boldsymbol{\dot{\Gamma}}\left(\boldsymbol{x}\right),\qquad\text{or} & \boldsymbol{\dot{\Gamma}}\left(\boldsymbol{x}\right)-\boldsymbol{\mathcal{\dot{P}}}^{T}\left(\boldsymbol{x}\right)\boldsymbol{\Gamma}\left(\boldsymbol{x}\right) & =\boldsymbol{\mathcal{P}}^{T}\left(\boldsymbol{x}\right)\boldsymbol{\dot{\Gamma}}\left(\boldsymbol{x}\right).\label{eq:Eq221}
\end{align}
Furthermore, from
\begin{align}
\boldsymbol{\mathcal{P}}^{T}\left(\boldsymbol{x}\right)\boldsymbol{\mathcal{Q}}^{T}\left(\boldsymbol{x}\right) & =\boldsymbol{0}
\end{align}
follows
\begin{align}
\boldsymbol{\mathcal{\dot{P}}}^{T}\left(\boldsymbol{x}\right)\boldsymbol{\mathcal{Q}}^{T}\left(\boldsymbol{x}\right) & =-\boldsymbol{\mathcal{P}}^{T}\left(\boldsymbol{x}\right)\boldsymbol{\mathcal{\dot{Q}}}^{T}\left(\boldsymbol{x}\right),\qquad\text{or} & \boldsymbol{\mathcal{\dot{P}}}^{T}\left(\boldsymbol{x}\right)\boldsymbol{\mathcal{Q}}^{T}\left(\boldsymbol{x}\right)\boldsymbol{\mathcal{Q}}^{T}\left(\boldsymbol{x}\right) & =-\boldsymbol{\mathcal{P}}^{T}\left(\boldsymbol{x}\right)\boldsymbol{\mathcal{\dot{Q}}}^{T}\left(\boldsymbol{x}\right)\boldsymbol{\mathcal{Q}}^{T}\left(\boldsymbol{x}\right)\label{eq:PdotQ}
\end{align}
due to the idempotence of projectors. Using Eqs. (\ref{eq:Eq221})
and (\ref{eq:PdotQ}) in Eq. (\ref{eq:PSLambdaDot_3}) yields 
\begin{align}
\varepsilon^{2}\boldsymbol{\Gamma}\left(\boldsymbol{x}\right)\boldsymbol{\ddot{x}} & =\varepsilon^{2}\boldsymbol{\mathcal{P}}^{T}\left(\boldsymbol{x}\right)\left(\boldsymbol{\Gamma}\left(\boldsymbol{x}\right)\nabla\boldsymbol{R}\left(\boldsymbol{x}\right)\boldsymbol{\dot{x}}-\boldsymbol{\dot{\Gamma}}\left(\boldsymbol{x}\right)\left(\boldsymbol{\dot{x}}-\boldsymbol{R}\left(\boldsymbol{x}\right)\right)\right)+\boldsymbol{\mathcal{P}}^{T}\left(\boldsymbol{x}\right)\boldsymbol{w}_{\varepsilon}\left(\boldsymbol{x},\boldsymbol{\dot{x}},\boldsymbol{\lambda}\right)\nonumber \\
 & +\boldsymbol{\mathcal{P}}^{T}\left(\boldsymbol{x}\right)\boldsymbol{\mathcal{\dot{Q}}}^{T}\left(\boldsymbol{x}\right)\boldsymbol{\mathcal{Q}}^{T}\left(\boldsymbol{x}\right)\boldsymbol{\lambda}+\boldsymbol{\mathcal{P}}^{T}\left(\boldsymbol{x}\right)\boldsymbol{\mathcal{S}}\left(\boldsymbol{x}-\boldsymbol{x}_{d}\left(t\right)\right).\label{eq:PSLambdaDot_4}
\end{align}
The form of Eq. (\ref{eq:PSLambdaDot_4}) together with 
\begin{align}
\boldsymbol{\mathcal{Q}}^{T}\left(\boldsymbol{x}\right)\boldsymbol{\Gamma}\left(\boldsymbol{x}\right)=\boldsymbol{\mathcal{Q}}^{T}\left(\boldsymbol{x}\right)\boldsymbol{\mathcal{P}}^{T}\left(\boldsymbol{x}\right)\boldsymbol{\Gamma}\left(\boldsymbol{x}\right) & =\boldsymbol{0}
\end{align}
makes it obvious that it contains no component in the ''direction''
$\boldsymbol{\mathcal{Q}}^{T}\left(\boldsymbol{x}\right)$. Equation
(\ref{eq:PSLambdaDot_4}) constitutes $p$ linearly independent second
order differential equations for $p$ linearly independent state components
$\boldsymbol{\mathcal{P}}\left(\boldsymbol{x}\right)\boldsymbol{x}$.
The $2p$ boundary conditions necessary to solve Eq. (\ref{eq:PSLambdaDot_4})
are given by Eqs. (\ref{eq:PSQSInitCond}) and (\ref{eq:PSQSTermCond}).\\
To summarize the derivation, the rearranged necessary optimality conditions
are 
\begin{align}
-\boldsymbol{\mathcal{Q}}^{T}\left(\boldsymbol{x}\right)\boldsymbol{\dot{\lambda}} & =\boldsymbol{\mathcal{Q}}^{T}\left(\boldsymbol{x}\right)\boldsymbol{w}_{\varepsilon}\left(\boldsymbol{x},\boldsymbol{\dot{x}},\boldsymbol{\lambda}\right)+\boldsymbol{\mathcal{Q}}^{T}\left(\boldsymbol{x}\right)\boldsymbol{\mathcal{S}}\boldsymbol{\mathcal{Q}}\left(\boldsymbol{x}\right)\left(\boldsymbol{x}-\boldsymbol{x}_{d}\left(t\right)\right),\label{eq:Rearranged1-1}\\
\varepsilon^{2}\boldsymbol{\Gamma}\left(\boldsymbol{x}\right)\boldsymbol{\ddot{x}} & =\varepsilon^{2}\boldsymbol{\mathcal{P}}^{T}\left(\boldsymbol{x}\right)\left(\boldsymbol{\Gamma}\left(\boldsymbol{x}\right)\nabla\boldsymbol{R}\left(\boldsymbol{x}\right)\boldsymbol{\dot{x}}-\boldsymbol{\dot{\Gamma}}\left(\boldsymbol{x}\right)\left(\boldsymbol{\dot{x}}-\boldsymbol{R}\left(\boldsymbol{x}\right)\right)\right)+\boldsymbol{\mathcal{P}}^{T}\left(\boldsymbol{x}\right)\boldsymbol{w}_{\varepsilon}\left(\boldsymbol{x},\boldsymbol{\dot{x}},\boldsymbol{\lambda}\right)\nonumber \\
 & +\boldsymbol{\mathcal{P}}^{T}\left(\boldsymbol{x}\right)\boldsymbol{\mathcal{\dot{Q}}}^{T}\left(\boldsymbol{x}\right)\boldsymbol{\mathcal{Q}}^{T}\left(\boldsymbol{x}\right)\boldsymbol{\lambda}+\boldsymbol{\mathcal{P}}^{T}\left(\boldsymbol{x}\right)\boldsymbol{\mathcal{S}}\boldsymbol{\mathcal{P}}\left(\boldsymbol{x}\right)\left(\boldsymbol{x}-\boldsymbol{x}_{d}\left(t\right)\right),\label{eq:Rearranged2-1}\\
\boldsymbol{\mathcal{Q}}\left(\boldsymbol{x}\right)\boldsymbol{\dot{x}} & =\boldsymbol{\mathcal{Q}}\left(\boldsymbol{x}\right)\boldsymbol{R}\left(\boldsymbol{x}\right),\label{eq:Rearranged3-1}
\end{align}
with $\boldsymbol{w}_{\varepsilon}$ defined in Eq. (\ref{eq:Defw}).
Equation (\ref{eq:SPCommute2}) was used for the terms $\boldsymbol{\mathcal{Q}}^{T}\boldsymbol{\mathcal{S}}\boldsymbol{\mathcal{Q}}$
and $\boldsymbol{\mathcal{P}}^{T}\boldsymbol{\mathcal{S}}\boldsymbol{\mathcal{P}}$.
We emphasize that Eqs. (\ref{eq:Rearranged1-1})-(\ref{eq:Rearranged3-1})
arise as a rearrangement of the necessary optimality conditions Eqs.
(\ref{eq:StationarityCondition-1})-(\ref{eq:AdjointEq-1}) and are
derived without approximation. The small regularization parameter
$\varepsilon$ multiplies the highest derivative $\boldsymbol{\ddot{x}}\left(t\right)$
in the system such that Eqs. (\ref{eq:Rearranged1-1})-(\ref{eq:Rearranged3-1})
constitute a system of singularly perturbed differential equations.

\section*{Appendix B: Outer equations and linearizing assumption}

The outer equations are obtained by setting $\varepsilon=0$ in Eqs.
(\ref{eq:Rearranged1-1})-(\ref{eq:Rearranged3-1}). Together with
the definition of $\boldsymbol{w}_{\varepsilon}$ in Eq. (\ref{eq:Defw}),
we obtain 
\begin{align}
-\boldsymbol{\mathcal{Q}}^{T}\left(\boldsymbol{x}\right)\boldsymbol{\dot{\lambda}} & =\boldsymbol{\mathcal{Q}}^{T}\left(\boldsymbol{x}\right)\left(\nabla\boldsymbol{R}^{T}\left(\boldsymbol{x}\right)+\boldsymbol{\mathcal{W}}^{T}\left(\boldsymbol{x},\boldsymbol{\dot{x}}\right)\right)\boldsymbol{\mathcal{Q}}^{T}\left(\boldsymbol{x}\right)\boldsymbol{\lambda}+\boldsymbol{\mathcal{Q}}^{T}\left(\boldsymbol{x}\right)\boldsymbol{\mathcal{S}}\boldsymbol{\mathcal{Q}}\left(\boldsymbol{x}\right)\left(\boldsymbol{x}-\boldsymbol{x}_{d}\left(t\right)\right),\label{eq:EqOut1}\\
-\boldsymbol{\mathcal{P}}^{T}\left(\boldsymbol{x}\right)\boldsymbol{\mathcal{\dot{Q}}}^{T}\left(\boldsymbol{x}\right)\boldsymbol{\mathcal{Q}}^{T}\left(\boldsymbol{x}\right)\boldsymbol{\lambda} & =\boldsymbol{\mathcal{P}}^{T}\left(\boldsymbol{x}\right)\left(\nabla\boldsymbol{R}^{T}\left(\boldsymbol{x}\right)+\boldsymbol{\mathcal{W}}^{T}\left(\boldsymbol{x},\boldsymbol{\dot{x}}\right)\right)\boldsymbol{\mathcal{Q}}^{T}\left(\boldsymbol{x}\right)\boldsymbol{\lambda}+\boldsymbol{\mathcal{P}}^{T}\left(\boldsymbol{x}\right)\boldsymbol{\mathcal{S}}\boldsymbol{\mathcal{P}}\left(\boldsymbol{x}\right)\left(\boldsymbol{x}-\boldsymbol{x}_{d}\left(t\right)\right),\label{eq:EqOut2}\\
\boldsymbol{\mathcal{Q}}\left(\boldsymbol{x}\right)\boldsymbol{\dot{x}} & =\boldsymbol{\mathcal{Q}}\left(\boldsymbol{x}\right)\boldsymbol{R}\left(\boldsymbol{x}\right).\label{eq:EqOut3}
\end{align}
Multiplying Eq. (\ref{eq:EqOut2}) by $\boldsymbol{\Omega}\left(\boldsymbol{x}\right)$
from the left and using Eqs. (\ref{eq:Eq4210}) and (\ref{eq:Eq4211})
yields
\begin{align}
\boldsymbol{\mathcal{P}}\left(\boldsymbol{x}\right)\boldsymbol{x} & =\boldsymbol{\mathcal{P}}\left(\boldsymbol{x}\right)\boldsymbol{x}_{d}\left(t\right)-\boldsymbol{\Omega}\left(\boldsymbol{x}\right)\left(\boldsymbol{\mathcal{\dot{Q}}}^{T}\left(\boldsymbol{x}\right)+\nabla\boldsymbol{R}^{T}\left(\boldsymbol{x}\right)+\boldsymbol{\mathcal{W}}^{T}\left(\boldsymbol{x},\boldsymbol{\dot{x}}\right)\right)\boldsymbol{\mathcal{Q}}^{T}\left(\boldsymbol{x}\right)\boldsymbol{\lambda}.\label{eq:EqOut21}
\end{align}
Before applying the linearizing assumption, we have to discuss how
the linearizing assumption affects the product $\boldsymbol{\mathcal{Q}}\left(\boldsymbol{x}\right)\boldsymbol{\mathcal{W}}\left(\boldsymbol{x},\boldsymbol{y}\right)$
with $\boldsymbol{\mathcal{W}}\left(\boldsymbol{x},\boldsymbol{y}\right)$
defined in Eq. (\ref{eq:DefW}). Using Eq. (\ref{eq:Eq4222}) yields
the identity
\begin{align}
\dfrac{\partial}{\partial x_{j}}\boldsymbol{\mathcal{B}}\left(\boldsymbol{x}\right)\boldsymbol{\mathcal{B}}^{g}\left(\boldsymbol{x}\right) & =\dfrac{\partial}{\partial x_{j}}\boldsymbol{\mathcal{P}}\left(\boldsymbol{x}\right)-\boldsymbol{\mathcal{B}}\left(\boldsymbol{x}\right)\dfrac{\partial}{\partial x_{j}}\boldsymbol{\mathcal{B}}^{g}\left(\boldsymbol{x}\right),
\end{align}
such that the entries of $\boldsymbol{\mathcal{W}}\left(\boldsymbol{x},\boldsymbol{y}\right)$
can be expressed as 
\begin{align}
\mathcal{W}_{ij}\left(\boldsymbol{x},\boldsymbol{y}\right) & =\sum_{k=1}^{n}\dfrac{\partial}{\partial x_{j}}\mathcal{P}_{ik}\left(\boldsymbol{x}\right)\left(y_{k}-R_{k}\left(\boldsymbol{x}\right)\right)-\sum_{k=1}^{n}\sum_{l=1}^{p}\mathcal{B}_{il}\left(\boldsymbol{x}\right)\dfrac{\partial}{\partial x_{j}}\mathcal{B}_{lk}^{g}\left(\boldsymbol{x}\right)\left(y_{k}-R_{k}\left(\boldsymbol{x}\right)\right).\label{eq:Wentries}
\end{align}
Because of $\boldsymbol{\mathcal{Q}}\left(\boldsymbol{x}\right)\boldsymbol{\mathcal{B}}\left(\boldsymbol{x}\right)=\boldsymbol{0}$,
the product $\boldsymbol{\mathcal{Q}}\left(\boldsymbol{x}\right)\boldsymbol{\mathcal{W}}\left(\boldsymbol{x},\boldsymbol{y}\right)$
is
\begin{align}
\sum_{i=1}^{n}\mathcal{Q}_{li}\left(\boldsymbol{x}\right)\mathcal{W}_{ij}\left(\boldsymbol{x},\boldsymbol{y}\right) & =\sum_{i=1}^{n}\sum_{k=1}^{n}\mathcal{Q}_{li}\left(\boldsymbol{x}\right)\dfrac{\partial}{\partial x_{j}}\mathcal{P}_{ik}\left(\boldsymbol{x}\right)\left(y_{k}-R_{k}\left(\boldsymbol{x}\right)\right).\label{eq:QTimesW}
\end{align}
Thus for a constant projector, $\boldsymbol{\mathcal{P}}\left(\boldsymbol{x}\right)=\boldsymbol{\mathcal{P}}=\text{const.}$,
the product $\boldsymbol{\mathcal{Q}}\left(\boldsymbol{x}\right)\boldsymbol{\mathcal{W}}\left(\boldsymbol{x},\boldsymbol{y}\right)=\boldsymbol{0}$
vanishes. For the sake of a clear notation, we denote the solution
to the outer equations by capital letters, $\boldsymbol{X}\left(t\right)=\boldsymbol{x}\left(t\right)$,
and $\boldsymbol{\Lambda}\left(t\right)=\boldsymbol{\lambda}\left(t\right)$.
Employing the linearizing assumption in Eqs. (\ref{eq:EqOut1}), (\ref{eq:EqOut3}),
and (\ref{eq:EqOut21}) finally yields the linear outer equations
\begin{align}
-\boldsymbol{\mathcal{Q}}^{T}\boldsymbol{\dot{\Lambda}}\left(t\right) & =\boldsymbol{\mathcal{Q}}^{T}\boldsymbol{\mathcal{A}}^{T}\boldsymbol{\mathcal{Q}}^{T}\boldsymbol{\Lambda}\left(t\right)+\boldsymbol{\mathcal{Q}}^{T}\boldsymbol{\mathcal{S}}\boldsymbol{\mathcal{Q}}\left(\boldsymbol{X}\left(t\right)-\boldsymbol{x}_{d}\left(t\right)\right),\label{eq:Eq91}\\
\boldsymbol{\mathcal{P}}\boldsymbol{X}\left(t\right) & =\boldsymbol{\mathcal{P}}\boldsymbol{x}_{d}\left(t\right)-\boldsymbol{\Omega}\boldsymbol{\mathcal{A}}^{T}\boldsymbol{\mathcal{Q}}^{T}\boldsymbol{\Lambda}\left(t\right),\label{eq:Eq92}\\
\boldsymbol{\mathcal{Q}}\boldsymbol{\dot{X}}\left(t\right) & =\boldsymbol{\mathcal{Q}}\boldsymbol{\mathcal{A}}\boldsymbol{X}\left(t\right)+\boldsymbol{\mathcal{Q}}\boldsymbol{b}.\label{eq:Eq93}
\end{align}
After elimination of $\boldsymbol{\mathcal{P}}\boldsymbol{X}\left(t\right)$
from Eq. (\ref{eq:Eq93}) with the help of Eq. (\ref{eq:Eq92}), we
can rewrite Eqs. (\ref{eq:Eq91}) and (\ref{eq:Eq93}) in matrix notation
as
\begin{align}
\left(\begin{array}{c}
\boldsymbol{\mathcal{Q}}^{T}\boldsymbol{\dot{\Lambda}}\left(t\right)\\
\boldsymbol{\mathcal{Q}}\boldsymbol{\dot{X}}\left(t\right)
\end{array}\right) & =\left(\begin{array}{cc}
-\boldsymbol{\mathcal{Q}}^{T}\boldsymbol{\mathcal{A}}^{T}\boldsymbol{\mathcal{Q}}^{T} & -\boldsymbol{\mathcal{Q}}^{T}\boldsymbol{\mathcal{S}}\boldsymbol{\mathcal{Q}}\\
-\boldsymbol{\mathcal{Q}}\boldsymbol{\mathcal{A}}\boldsymbol{\mathcal{P}}\boldsymbol{\Omega}\boldsymbol{\mathcal{A}}^{T}\boldsymbol{\mathcal{Q}}^{T} & \boldsymbol{\mathcal{Q}}\boldsymbol{\mathcal{A}}\boldsymbol{\mathcal{Q}}
\end{array}\right)\left(\begin{array}{c}
\boldsymbol{\mathcal{Q}}^{T}\boldsymbol{\Lambda}\left(t\right)\\
\boldsymbol{\mathcal{Q}}\boldsymbol{X}\left(t\right)
\end{array}\right)+\left(\begin{array}{c}
\boldsymbol{\mathcal{Q}}^{T}\boldsymbol{\mathcal{S}}\boldsymbol{\mathcal{Q}}\boldsymbol{x}_{d}\left(t\right)\\
\boldsymbol{\mathcal{Q}}\boldsymbol{\mathcal{A}}\boldsymbol{\mathcal{P}}\boldsymbol{x}_{d}\left(t\right)+\boldsymbol{\mathcal{Q}}\boldsymbol{b}
\end{array}\right).\label{eq:OuterEquationsS}
\end{align}
The are $2\left(n-p\right)$ linearly independent ODEs for the $2\left(n-p\right)$
linearly independent components of $\boldsymbol{\mathcal{Q}}^{T}\boldsymbol{\Lambda}\left(t\right)$
and $\boldsymbol{\mathcal{Q}}\boldsymbol{X}\left(t\right)$.

\section*{\label{sec:AppendixB}Appendix C: Inner equations}

The left inner equations, valid near to the initial time $t\gtrsim t_{0}$,
are resolved by the time scale $\tau_{L}$ defined as
\begin{align}
\tau_{L} & =\left(t-t_{0}\right)/\varepsilon.
\end{align}
The left inner solutions are denoted by capital letters with index
$L$, 
\begin{align}
\boldsymbol{X}_{L}\left(\tau_{L}\right) & =\boldsymbol{X}_{L}\left(\left(t-t_{0}\right)/\varepsilon\right)=\boldsymbol{x}\left(t\right), & \boldsymbol{\Lambda}_{L}\left(\tau_{L}\right) & =\boldsymbol{\Lambda}_{L}\left(\left(t-t_{0}\right)/\varepsilon\right)=\boldsymbol{\lambda}\left(t\right).
\end{align}
Expressed in terms of the inner solutions, the time derivatives become
\begin{align}
\boldsymbol{\dot{x}}\left(t\right) & =\frac{d}{dt}\boldsymbol{X}_{L}\left(\left(t-t_{0}\right)/\varepsilon\right)=\varepsilon^{-1}\boldsymbol{X}_{L}'\left(\tau_{L}\right), & \boldsymbol{\dot{\lambda}}\left(t\right) & =\varepsilon^{-1}\boldsymbol{\Lambda}_{L}'\left(\tau_{L}\right),\\
\boldsymbol{\ddot{x}}\left(t\right) & =\frac{d^{2}}{dt^{2}}\boldsymbol{X}_{L}\left(\left(t-t_{0}\right)/\varepsilon\right)=\varepsilon^{-2}\boldsymbol{X}_{L}''\left(\tau_{L}\right).
\end{align}
The prime $\left(\cdot\right)'$ denotes the derivative with respect
to $\tau_{L}$, $\left(\cdot\right)'=\dfrac{\partial}{\partial\tau_{L}}$$\left(\cdot\right)$.
The time derivatives of $\boldsymbol{\mathcal{Q}}\left(\boldsymbol{x}\left(t\right)\right)$
and $\boldsymbol{\Gamma}\left(\boldsymbol{x}\left(t\right)\right)$
transform as 
\begin{align}
\boldsymbol{\mathcal{\dot{Q}}}^{T}\left(\boldsymbol{x}\left(t\right)\right) & =\varepsilon^{-1}\boldsymbol{\mathcal{Q}}^{T}\vspace{0cm}'\left(\boldsymbol{X}_{L}\left(\tau_{L}\right)\right), & \boldsymbol{\dot{\Gamma}}\left(\boldsymbol{x}\left(t\right)\right) & =\varepsilon^{-1}\boldsymbol{\Gamma}'\left(\boldsymbol{X}_{L}\left(\tau_{L}\right)\right).
\end{align}
The matrix $\boldsymbol{\mathcal{W}}\left(\boldsymbol{x}\left(t\right),\boldsymbol{\dot{x}}\left(t\right)\right)$
transforms as 
\begin{align}
\boldsymbol{\mathcal{W}}\left(\boldsymbol{x}\left(t\right),\boldsymbol{\dot{x}}\left(t\right)\right) & =\nabla\boldsymbol{\mathcal{B}}\left(\boldsymbol{x}\left(t\right)\right)\boldsymbol{\mathcal{B}}^{g}\left(\boldsymbol{x}\left(t\right)\right)\boldsymbol{\dot{x}}\left(t\right)-\nabla\boldsymbol{\mathcal{B}}\left(\boldsymbol{x}\left(t\right)\right)\boldsymbol{\mathcal{B}}^{g}\left(\boldsymbol{x}\left(t\right)\right)\boldsymbol{R}\left(\boldsymbol{x}\left(t\right)\right)\nonumber \\
 & =\varepsilon^{-1}\boldsymbol{\mathcal{V}}\left(\boldsymbol{X}_{L}\left(\tau_{L}\right),\boldsymbol{X}_{L}'\left(\tau_{L}\right)\right)+\boldsymbol{\mathcal{U}}\left(\boldsymbol{X}_{L}\left(\tau_{L}\right)\right),
\end{align}
with $n\times n$ matrices $\boldsymbol{\mathcal{U}}$ and $\boldsymbol{\mathcal{V}}$
defined by 
\begin{align}
\boldsymbol{\mathcal{U}}\left(\boldsymbol{x}\right) & =-\nabla\boldsymbol{\mathcal{B}}\left(\boldsymbol{x}\right)\boldsymbol{\mathcal{B}}^{g}\left(\boldsymbol{x}\right)\boldsymbol{R}\left(\boldsymbol{x}\right), & \boldsymbol{\mathcal{V}}\left(\boldsymbol{x},\boldsymbol{y}\right) & =\nabla\boldsymbol{\mathcal{B}}\left(\boldsymbol{x}\right)\boldsymbol{\mathcal{B}}^{g}\left(\boldsymbol{x}\right)\boldsymbol{y}.
\end{align}
The entries of $\boldsymbol{\mathcal{U}}$ and $\boldsymbol{\mathcal{V}}$
are 
\begin{align}
\mathcal{U}_{ij}\left(\boldsymbol{x}\right) & =\sum_{k=1}^{n}\sum_{l=1}^{p}\dfrac{\partial}{\partial x_{j}}\mathcal{B}_{il}\left(\boldsymbol{x}\right)\mathcal{B}_{lk}^{g}\left(\boldsymbol{x}\right)R_{k}\left(\boldsymbol{x}\right), & \mathcal{V}_{ij}\left(\boldsymbol{x},\boldsymbol{y}\right) & =\sum_{k=1}^{n}\sum_{l=1}^{p}\dfrac{\partial}{\partial x_{j}}\mathcal{B}_{il}\left(\boldsymbol{x}\right)\mathcal{B}_{lk}^{g}\left(\boldsymbol{x}\right)y_{k}.
\end{align}
From the initial conditions Eq. (\ref{eq:GenDynSysInitTermCond})
follow the initial conditions for $\boldsymbol{X}_{L}\left(\tau_{L}\right)$
as
\begin{align}
\boldsymbol{X}_{L}\left(0\right) & =\boldsymbol{x}_{0}.\label{eq:InitCondL}
\end{align}
Transforming the necessary optimality conditions Eqs. (\ref{eq:Rearranged1-1})-(\ref{eq:Rearranged3-1})
yields
\begin{align}
-\varepsilon^{-1}\boldsymbol{\mathcal{Q}}^{T}\left(\boldsymbol{X}_{L}\right)\boldsymbol{\Lambda}_{L}' & =\boldsymbol{\mathcal{Q}}^{T}\left(\boldsymbol{X}_{L}\right)\boldsymbol{\mathcal{V}}^{T}\left(\boldsymbol{X}_{L},\boldsymbol{X}_{L}'\right)\left(\boldsymbol{\Gamma}\left(\boldsymbol{X}_{L}\right)\left(\varepsilon\boldsymbol{R}\left(\boldsymbol{X}_{L}\right)-\boldsymbol{X}_{L}'\right)+\varepsilon^{-1}\boldsymbol{\mathcal{Q}}^{T}\left(\boldsymbol{X}_{L}\right)\boldsymbol{\Lambda}_{L}\right)\nonumber \\
 & +\boldsymbol{\mathcal{Q}}^{T}\left(\boldsymbol{X}_{L}\right)\left(\nabla\boldsymbol{R}^{T}\left(\boldsymbol{X}_{L}\right)+\boldsymbol{\mathcal{U}}^{T}\left(\boldsymbol{X}_{L}\right)\right)\left(\boldsymbol{\Gamma}\left(\boldsymbol{X}_{L}\right)\left(\varepsilon^{2}\boldsymbol{R}\left(\boldsymbol{X}_{L}\right)-\varepsilon\boldsymbol{X}_{L}'\right)+\boldsymbol{\mathcal{Q}}^{T}\left(\boldsymbol{X}_{L}\right)\boldsymbol{\Lambda}_{L}\right)\nonumber \\
 & +\boldsymbol{\mathcal{Q}}^{T}\left(\boldsymbol{X}_{L}\right)\boldsymbol{\mathcal{S}}\boldsymbol{\mathcal{Q}}\left(\boldsymbol{X}_{L}\right)\left(\boldsymbol{X}_{L}-\boldsymbol{x}_{d}\left(t_{0}+\varepsilon\tau_{L}\right)\right),\label{eq:InnerLEq1General}\\
\boldsymbol{\Gamma}\left(\boldsymbol{X}_{L}\right)\boldsymbol{X}_{L}'' & =\boldsymbol{\mathcal{P}}^{T}\left(\boldsymbol{X}_{L}\right)\left(\boldsymbol{\mathcal{V}}^{T}\left(\boldsymbol{X}_{L},\boldsymbol{X}_{L}'\right)\boldsymbol{\Gamma}\left(\boldsymbol{X}_{L}\right)+\boldsymbol{\Gamma}'\left(\boldsymbol{X}_{L}\right)\right)\left(\varepsilon\boldsymbol{R}\left(\boldsymbol{X}_{L}\right)-\boldsymbol{X}_{L}'\right)\nonumber \\
 & +\varepsilon^{-1}\boldsymbol{\mathcal{P}}^{T}\left(\boldsymbol{X}_{L}\right)\left(\boldsymbol{\mathcal{V}}^{T}\left(\boldsymbol{X}_{L},\boldsymbol{X}_{L}'\right)+\boldsymbol{\mathcal{Q}}^{T}\vspace{0cm}'\left(\boldsymbol{X}_{L}\right)\right)\boldsymbol{\mathcal{Q}}^{T}\left(\boldsymbol{X}_{L}\right)\boldsymbol{\Lambda}_{L}\nonumber \\
 & +\boldsymbol{\mathcal{P}}^{T}\left(\boldsymbol{X}_{L}\right)\left(\nabla\boldsymbol{R}^{T}\left(\boldsymbol{X}_{L}\right)+\boldsymbol{\mathcal{U}}^{T}\left(\boldsymbol{X}_{L}\right)\right)\left(\boldsymbol{\Gamma}\left(\boldsymbol{X}_{L}\right)\left(\varepsilon^{2}\boldsymbol{R}\left(\boldsymbol{X}_{L}\right)-\varepsilon\boldsymbol{X}_{L}'\right)+\boldsymbol{\mathcal{Q}}^{T}\left(\boldsymbol{X}_{L}\right)\boldsymbol{\Lambda}_{L}\right)\nonumber \\
 & -\varepsilon\boldsymbol{\mathcal{P}}^{T}\left(\boldsymbol{X}_{L}\right)\boldsymbol{\Gamma}\left(\boldsymbol{X}_{L}\right)\nabla\boldsymbol{R}\left(\boldsymbol{X}_{L}\right)\boldsymbol{X}_{L}'\left(\tau_{L}\right)+\boldsymbol{\mathcal{P}}^{T}\left(\boldsymbol{X}_{L}\right)\boldsymbol{\mathcal{S}}\boldsymbol{\mathcal{P}}\left(\boldsymbol{X}_{L}\right)\left(\boldsymbol{X}_{L}-\boldsymbol{x}_{d}\left(t_{0}+\varepsilon\tau_{L}\right)\right),\label{eq:InnerLEq2General}\\
\varepsilon^{-1}\boldsymbol{\mathcal{Q}}\left(\boldsymbol{X}_{L}\right)\boldsymbol{X}_{L}' & =\boldsymbol{\mathcal{Q}}\left(\boldsymbol{X}_{L}\right)\boldsymbol{R}\left(\boldsymbol{X}_{L}\right).\label{eq:InnerLEq3General}
\end{align}
In the next next step we enforce the linearizing assumption and set
$\varepsilon=0$. The linearizing assumption implies constant projectors,
$\boldsymbol{\mathcal{P}}\left(\boldsymbol{X}_{L}\right)=\boldsymbol{\mathcal{P}}$
and $\boldsymbol{\mathcal{Q}}\left(\boldsymbol{X}_{L}\right)=\boldsymbol{\mathcal{Q}}$,
and 
\begin{align}
\boldsymbol{\mathcal{V}}^{T}\left(\boldsymbol{X}_{L},\boldsymbol{X}_{L}'\right)\boldsymbol{\mathcal{Q}}^{T} & =\mathbf{0}, & \boldsymbol{\mathcal{U}}^{T}\left(\boldsymbol{X}_{L}\right)\boldsymbol{\mathcal{Q}}^{T} & =\mathbf{0},
\end{align}
and Eqs. (\ref{eq:InnerLEq1General})-(\ref{eq:InnerLEq3General})
become 
\begin{align}
\boldsymbol{\mathcal{Q}}^{T}\boldsymbol{\Lambda}_{L}' & =\mathbf{0},\\
\dfrac{\partial}{\partial\tau_{L}}\left(\boldsymbol{\Gamma}\left(\boldsymbol{X}_{L}\right)\boldsymbol{X}_{L}'\right) & =\boldsymbol{\mathcal{P}}^{T}\boldsymbol{\mathcal{A}}^{T}\boldsymbol{\mathcal{Q}}^{T}\boldsymbol{\Lambda}_{L}-\boldsymbol{\mathcal{P}}^{T}\boldsymbol{\mathcal{V}}^{T}\left(\boldsymbol{X}_{L},\boldsymbol{X}_{L}'\right)\boldsymbol{\Gamma}\left(\boldsymbol{X}_{L}\right)\boldsymbol{X}_{L}'+\boldsymbol{\mathcal{P}}^{T}\boldsymbol{\mathcal{S}}\boldsymbol{\mathcal{P}}\left(\boldsymbol{X}_{L}-\boldsymbol{x}_{d}\left(t_{0}\right)\right),\\
\boldsymbol{\mathcal{Q}}\boldsymbol{X}_{L}' & =\mathbf{0}.\label{eq:QSXLSol}
\end{align}
The right inner equations, valid near to the terminal time $t\lesssim t_{1}$,
are similarly dealt with as the left inner equations. The new time
scale is
\begin{align}
\tau_{R} & =\left(t_{1}-t\right)/\varepsilon.
\end{align}
The inner solutions are denoted by capital letters with an index $R$,
\begin{align}
\boldsymbol{X}_{R}\left(\tau_{R}\right) & =\boldsymbol{X}_{R}\left(\left(t_{1}-t\right)/\varepsilon\right)=\boldsymbol{x}\left(t\right), & \boldsymbol{\Lambda}_{R}\left(\tau_{R}\right) & =\boldsymbol{\Lambda}_{R}\left(\left(t_{1}-t\right)/\varepsilon\right)=\boldsymbol{\lambda}\left(t\right).
\end{align}
The derivation of the leading order right inner equations proceeds
analogous to the inner equations on the left side. The only difference
is that a minus sign appears for time derivatives of odd order. Note
that $\boldsymbol{\mathcal{V}}\left(\boldsymbol{x},-\boldsymbol{y}\right)=-\boldsymbol{\mathcal{V}}\left(\boldsymbol{x},\boldsymbol{y}\right)$.
Together with the linearizing assumption, we obtain
\begin{align}
\boldsymbol{\mathcal{Q}}^{T}\boldsymbol{\Lambda}_{R}' & =-\boldsymbol{\mathcal{Q}}^{T}\boldsymbol{\mathcal{V}}^{T}\left(\boldsymbol{X}_{R},\boldsymbol{X}_{R}'\right)\boldsymbol{\mathcal{Q}}^{T}\boldsymbol{\Lambda}_{R},\\
\dfrac{\partial}{\partial\tau_{R}}\left(\boldsymbol{\Gamma}\left(\boldsymbol{X}_{R}\right)\boldsymbol{X}_{R}'\right) & =\boldsymbol{\mathcal{P}}^{T}\boldsymbol{\mathcal{A}}^{T}\boldsymbol{\mathcal{Q}}^{T}\boldsymbol{\Lambda}_{R}-\boldsymbol{\mathcal{P}}^{T}\boldsymbol{\mathcal{V}}^{T}\left(\boldsymbol{X}_{R},\boldsymbol{X}_{R}'\right)\boldsymbol{\Gamma}\left(\boldsymbol{X}_{R}\right)\boldsymbol{X}_{R}'+\boldsymbol{\mathcal{P}}^{T}\boldsymbol{\mathcal{S}}\boldsymbol{\mathcal{P}}\left(\boldsymbol{X}_{R}-\boldsymbol{x}_{d}\left(t_{1}\right)\right),\\
\boldsymbol{\mathcal{Q}}\boldsymbol{X}_{R}' & =\boldsymbol{0}.\label{eq:QSXRSol}
\end{align}

\section*{\label{sec:AppendixC}Appendix D: Matching, composite solution, and
exact solution for $\varepsilon=0$}

The left and right inner equations and the outer equations must be
solved with appropriate boundary and matching conditions \cite{bender1999advanced}.
The boundary conditions for the state, Eq. (\ref{eq:GenDynSysInitTermCond}),
transform to boundary conditions for the inner equations, 
\begin{align}
\boldsymbol{\mathcal{P}}\boldsymbol{X}_{L}\left(0\right) & =\boldsymbol{\mathcal{P}}\boldsymbol{x}_{0}, & \boldsymbol{\mathcal{Q}}\boldsymbol{X}_{L}\left(0\right) & =\boldsymbol{\mathcal{Q}}\boldsymbol{x}_{0}, & \boldsymbol{\mathcal{P}}\boldsymbol{X}_{R}\left(0\right) & =\boldsymbol{\mathcal{P}}\boldsymbol{x}_{1}, & \boldsymbol{\mathcal{Q}}\boldsymbol{X}_{R}\left(0\right) & =\boldsymbol{\mathcal{Q}}\boldsymbol{x}_{1}.
\end{align}
On the left side, the matching conditions are
\begin{align}
\lim_{\tau_{L}\rightarrow\infty}\boldsymbol{\mathcal{Q}}^{T}\boldsymbol{\Lambda}_{L}\left(\tau_{L}\right) & =\lim_{t\rightarrow t{}_{0}}\boldsymbol{\mathcal{Q}}^{T}\boldsymbol{\Lambda}\left(t\right), & \lim_{\tau_{L}\rightarrow\infty}\boldsymbol{\mathcal{Q}}\boldsymbol{X}_{L}\left(\tau_{L}\right) & =\lim_{t\rightarrow t{}_{0}}\boldsymbol{\mathcal{Q}}\boldsymbol{X}\left(t\right), & \lim_{\tau_{L}\rightarrow\infty}\boldsymbol{\mathcal{P}}\boldsymbol{X}_{L}\left(\tau_{L}\right) & =\lim_{t\rightarrow t{}_{0}}\boldsymbol{\mathcal{P}}\boldsymbol{X}\left(t\right).\label{eq:LeftMatching}
\end{align}
Analogously, the matching conditions on the right side are
\begin{align}
\lim_{\tau_{R}\rightarrow\infty}\boldsymbol{\mathcal{Q}}^{T}\boldsymbol{\Lambda}_{R}\left(\tau_{R}\right) & =\lim_{t\rightarrow t{}_{1}}\boldsymbol{\mathcal{Q}}^{T}\boldsymbol{\Lambda}\left(t\right), & \lim_{\tau_{R}\rightarrow\infty}\boldsymbol{\mathcal{Q}}\boldsymbol{X}_{R}\left(\tau_{R}\right) & =\lim_{t\rightarrow t{}_{1}}\boldsymbol{\mathcal{Q}}\boldsymbol{X}\left(t\right), & \lim_{\tau_{R}\rightarrow\infty}\boldsymbol{\mathcal{P}}\boldsymbol{X}_{R}\left(\tau_{R}\right) & =\lim_{t\rightarrow t{}_{1}}\boldsymbol{\mathcal{P}}\boldsymbol{X}\left(t\right).\label{eq:RightMatching}
\end{align}
Evaluating the algebraic outer equation, Eq. (\ref{eq:Eq92}), at
the initial and terminal times $t_{0}$ and $t_{1}$ together with
the left and right matching conditions, Eqs. (\ref{eq:LeftMatching})
and (\ref{eq:RightMatching}), results in boundary conditions for
$\boldsymbol{\mathcal{P}}\boldsymbol{X}_{L}$ and $\boldsymbol{\mathcal{P}}\boldsymbol{X}_{R}$,
\begin{align}
\lim_{\tau_{L}\rightarrow\infty}\boldsymbol{\mathcal{P}}\boldsymbol{X}_{L}\left(\tau_{L}\right) & =\boldsymbol{\mathcal{P}}\boldsymbol{x}_{d}\left(t_{0}\right)-\boldsymbol{\Omega}\boldsymbol{\mathcal{A}}^{T}\boldsymbol{\mathcal{Q}}^{T}\boldsymbol{\Lambda}\left(t_{0}\right), & \lim_{\tau_{R}\rightarrow\infty}\boldsymbol{\mathcal{P}}\boldsymbol{X}_{R}\left(\tau_{R}\right) & =\boldsymbol{\mathcal{P}}\boldsymbol{x}_{d}\left(t_{1}\right)-\boldsymbol{\Omega}\boldsymbol{\mathcal{A}}^{T}\boldsymbol{\mathcal{Q}}^{T}\boldsymbol{\Lambda}\left(t_{1}\right).\label{eq:InnerLTermCond}
\end{align}
The existence of these limits, together with the result that $\boldsymbol{\mathcal{Q}}\boldsymbol{X}_{L}$
and $\boldsymbol{\mathcal{Q}}\boldsymbol{X}_{R}$ is constant, see
Eqs. (\ref{eq:QSXLSol}) and (\ref{eq:QSXRSol}), implies 
\begin{align}
\lim_{\tau_{L}\rightarrow\infty}\boldsymbol{X}_{L}'\left(\tau_{L}\right) & =\boldsymbol{0}, & \lim_{\tau_{R}\rightarrow\infty}\boldsymbol{X}_{R}'\left(\tau_{R}\right) & =\boldsymbol{0}.\label{eq:XLRTauLInfty}
\end{align}
Finally, the two matching conditions
\begin{align}
\lim_{\tau_{L}\rightarrow\infty}\boldsymbol{\mathcal{Q}}\boldsymbol{X}_{L}\left(\tau_{L}\right) & =\boldsymbol{\mathcal{Q}}\boldsymbol{X}\left(t_{0}\right), & \lim_{\tau_{R}\rightarrow\infty}\boldsymbol{\mathcal{Q}}\boldsymbol{X}_{R}\left(\tau_{R}\right) & =\boldsymbol{\mathcal{Q}}\boldsymbol{X}\left(t_{1}\right),
\end{align}
remain. Because of the constancy of $\boldsymbol{\mathcal{Q}}\boldsymbol{X}_{L}$
and $\boldsymbol{\mathcal{Q}}\boldsymbol{X}_{R}$, Eqs. (\ref{eq:QSXLSol})
and (\ref{eq:QSXRSol}), together with the matching conditions, Eqs.
(\ref{eq:LeftMatching}) and (\ref{eq:RightMatching}), respectively,
these can be written as 
\begin{align}
\boldsymbol{\mathcal{Q}}\boldsymbol{X}\left(t_{0}\right) & =\boldsymbol{\mathcal{Q}}\boldsymbol{x}_{0}, & \boldsymbol{\mathcal{Q}}\boldsymbol{X}\left(t_{1}\right) & =\boldsymbol{\mathcal{Q}}\boldsymbol{x}_{1}.
\end{align}
These are the $2\left(n-p\right)$ linearly independent boundary conditions
for the outer equations, Eq. (\ref{eq:OuterEquationsS}). They depend
only on the initial and terminal conditions $\boldsymbol{x}_{0}$
and $\boldsymbol{x}_{1}$, respectively. Hence, the solutions to the
outer equations are independent of any details of the inner equations.\\
Eventually, it is possible to formally write down the composite solutions
$\boldsymbol{x}_{\text{comp}}\left(t\right)$ and $\boldsymbol{\lambda}_{\text{comp}}\left(t\right)$
for the problem. The parts $\boldsymbol{\mathcal{Q}}\boldsymbol{x}_{\text{comp}}\left(t\right)$
and $\boldsymbol{\mathcal{Q}}\boldsymbol{\lambda}_{\text{comp}}\left(t\right)$
do not exhibit boundary layers and are simply given by the solutions
to the outer equations,
\begin{align}
\boldsymbol{\mathcal{Q}}\boldsymbol{x}_{\text{comp}}\left(t\right) & =\boldsymbol{\mathcal{Q}}\boldsymbol{X}\left(t\right), & \boldsymbol{\mathcal{Q}}\boldsymbol{\lambda}_{\text{comp}}\left(t\right) & =\boldsymbol{\mathcal{Q}}\boldsymbol{\Lambda}\left(t\right).
\end{align}
The part $\boldsymbol{\mathcal{P}}\boldsymbol{x}_{\text{comp}}\left(t\right)$
contains boundary layers and is given by the sum of outer, left inner
and right inner solution minus the overlaps $\boldsymbol{\mathcal{P}}\boldsymbol{X}\left(t_{0}\right)$
and $\boldsymbol{\mathcal{P}}\boldsymbol{X}\left(t_{1}\right)$ \cite{bender1999advanced},
\begin{align}
\boldsymbol{\mathcal{P}}\boldsymbol{x}_{\text{comp}}\left(t\right) & =\boldsymbol{\mathcal{P}}\boldsymbol{X}\left(t\right)+\boldsymbol{\mathcal{P}}\boldsymbol{X}_{L}\left(\left(t-t_{0}\right)/\varepsilon\right)+\boldsymbol{\mathcal{P}}\boldsymbol{X}_{R}\left(\left(t_{1}-t\right)/\varepsilon\right)-\boldsymbol{\mathcal{P}}\boldsymbol{X}\left(t_{0}\right)-\boldsymbol{\mathcal{P}}\boldsymbol{X}\left(t_{1}\right).\label{eq:CompSolGen3}
\end{align}
The controlled state reads as
\begin{align}
\boldsymbol{x}_{\text{comp}}\left(t\right) & =\boldsymbol{\mathcal{P}}\boldsymbol{x}_{\text{comp}}\left(t\right)+\boldsymbol{\mathcal{Q}}\boldsymbol{x}_{\text{comp}}\left(t\right)\nonumber \\
 & =\boldsymbol{X}\left(t\right)-\boldsymbol{\mathcal{P}}\boldsymbol{X}\left(t_{0}\right)-\boldsymbol{\mathcal{P}}\boldsymbol{X}\left(t_{1}\right)+\boldsymbol{\mathcal{P}}\boldsymbol{X}_{L}\left(\left(t-t_{0}\right)/\varepsilon\right)+\boldsymbol{\mathcal{P}}\boldsymbol{X}_{R}\left(\left(t_{1}-t\right)/\varepsilon\right).
\end{align}
The composite control signal is given in terms of the composite solutions
as
\begin{align}
\boldsymbol{u}_{\text{comp}}\left(t\right) & =\boldsymbol{\mathcal{B}}^{g}\left(\boldsymbol{x}_{\text{comp}}\left(t\right)\right)\left(\boldsymbol{\dot{x}}_{\text{comp}}\left(t\right)-\boldsymbol{R}\left(\boldsymbol{x}_{\text{comp}}\left(t\right)\right)\right).\label{eq:CompositeControl}
\end{align}
The exact state solution for $\varepsilon=0$ is obtained by taking
the limit $\varepsilon\rightarrow0$ of the composite solutions $\boldsymbol{x}_{\text{comp}}\left(t\right)$
and $\boldsymbol{\lambda}_{\text{comp}}\left(t\right)$. Only $\boldsymbol{\mathcal{P}}\boldsymbol{x}_{\text{comp}}\left(t\right)$
depends on $\varepsilon$, which yields
\begin{align}
\lim_{\varepsilon\rightarrow0}\boldsymbol{\mathcal{P}}\boldsymbol{x}_{\text{comp}}\left(t\right) & =\begin{cases}
\boldsymbol{\mathcal{P}}\boldsymbol{x}_{0}, & t=t_{0},\\
\boldsymbol{\mathcal{P}}\boldsymbol{X}\left(t\right), & t_{0}<t<t_{1},\\
\boldsymbol{\mathcal{P}}\boldsymbol{x}_{1}, & t=t_{1}.
\end{cases}
\end{align}
To obtain the exact control solution for $\varepsilon=0$, we have
to analyze Eq. (\ref{eq:CompositeControl}) in the limit $\varepsilon\rightarrow0$.
All terms except $\boldsymbol{\dot{x}}_{\text{comp}}\left(t\right)$
are well behaved. The term $\boldsymbol{\dot{x}}_{\text{comp}}\left(t\right)$
requires the investigation of the limit
\begin{align}
\lim_{\varepsilon\rightarrow0}\boldsymbol{\dot{x}}_{\text{comp}}\left(t\right) & =\boldsymbol{\dot{X}}\left(t\right)+\lim_{\varepsilon\rightarrow0}\boldsymbol{\mathcal{P}}\dfrac{d}{dt}\boldsymbol{X}_{L}\left(\left(t-t_{0}\right)/\varepsilon\right)+\lim_{\varepsilon\rightarrow0}\boldsymbol{\mathcal{P}}\dfrac{d}{dt}\boldsymbol{X}_{R}\left(\left(t_{1}-t\right)/\varepsilon\right).
\end{align}
The first step is to prove that $\boldsymbol{\mathcal{P}}\dfrac{d}{dt}\boldsymbol{X}_{L}\left(t/\varepsilon\right)$
yields a term proportional to the Dirac delta function $\delta\left(t\right)$
in the limit $\varepsilon\rightarrow0$. Define the $n\times1$ vector
of functions with 
\begin{align}
\boldsymbol{\delta}_{L,\varepsilon}\left(t\right) & =\begin{cases}
\boldsymbol{\mathcal{P}}\dfrac{d}{dt}\boldsymbol{X}_{L}\left(t/\varepsilon\right), & t\geq0,\\
\boldsymbol{\mathcal{P}}\dfrac{d}{d\tilde{t}}\boldsymbol{X}_{L}\left(\tilde{t}/\varepsilon\right)\Bigg|_{_{\tilde{t}=-t}}, & t<0.
\end{cases}
\end{align}
The function $\boldsymbol{\delta}_{L,\varepsilon}\left(t\right)$
is continuous for $t=0$ in every component. It can also be expressed
as
\begin{align}
\boldsymbol{\delta}_{L,\varepsilon}\left(t\right) & =\varepsilon^{-1}\boldsymbol{\mathcal{P}}\boldsymbol{X}_{L}'\left(\left|t\right|/\varepsilon\right).\label{eq:DefDelta}
\end{align}
First, evaluating $\boldsymbol{\delta}_{L,\varepsilon}\left(t\right)$
at $t=0$ yields 
\begin{align}
\boldsymbol{\delta}_{L,\varepsilon}\left(0\right) & =\varepsilon^{-1}\boldsymbol{\mathcal{P}}\boldsymbol{X}_{L}'\left(0\right),
\end{align}
and because $\boldsymbol{X}_{L}'\left(0\right)$ is finite and does
not depend on $\varepsilon$, this expression clearly diverges in
the limit $\varepsilon\rightarrow0$. Second, for $\left|t\right|>0$,
$\lim_{\varepsilon\rightarrow0}\boldsymbol{\delta}_{L,\varepsilon}\left(t\right)$
behaves as 
\begin{align}
\lim_{\varepsilon\rightarrow0}\boldsymbol{\delta}_{L,\varepsilon}\left(t\right) & =\boldsymbol{0},\,t\neq0,
\end{align}
because $\boldsymbol{X}_{L}'\left(\left|t\right|/\varepsilon\right)$
behaves as (see also Eq. (\ref{eq:XLRTauLInfty}))
\begin{align}
\lim_{\varepsilon\rightarrow0}\boldsymbol{X}_{L}'\left(\left|t\right|/\varepsilon\right) & =0,\,t\neq0.
\end{align}
Third, the integral of $\boldsymbol{\delta}_{L,\varepsilon}\left(t\right)$
over time $t$ must be determined. The integral can be split up in
two integrals,
\begin{align}
\intop_{-\infty}^{\infty}d\tilde{t}\boldsymbol{\delta}_{L,\varepsilon}\left(\tilde{t}\right) & =\intop_{-\infty}^{0}d\tilde{t}\boldsymbol{\delta}_{L,\varepsilon}\left(\tilde{t}\right)+\intop_{0}^{\infty}d\tilde{t}\boldsymbol{\delta}_{L,\varepsilon}\left(\tilde{t}\right)=\varepsilon^{-1}\intop_{-\infty}^{0}d\tilde{t}\boldsymbol{\mathcal{P}}\boldsymbol{X}_{L}'\left(-\tilde{t}/\varepsilon\right)+\varepsilon^{-1}\intop_{0}^{\infty}d\tilde{t}\boldsymbol{\mathcal{P}}\boldsymbol{X}_{L}'\left(\tilde{t}/\varepsilon\right).
\end{align}
Substituting $\tau=-\tilde{t}/\varepsilon$ in the first and $\tau=\tilde{t}/\varepsilon$
in the second integral yields (see Eq. (\ref{eq:LeftMatching})) 
\begin{align}
\intop_{-\infty}^{\infty}d\tilde{t}\boldsymbol{\delta}_{L,\varepsilon}\left(\tilde{t}\right) & =2\intop_{0}^{\infty}d\tau\boldsymbol{\mathcal{P}}\boldsymbol{X}_{L}'\left(\tau\right)=2\boldsymbol{\mathcal{P}}\left(\boldsymbol{X}\left(t_{0}\right)-\boldsymbol{x}_{0}\right).
\end{align}
Thus, we proved that
\begin{align}
\lim_{\varepsilon\rightarrow0}\boldsymbol{\delta}_{L,\varepsilon}\left(t\right) & =2\boldsymbol{\mathcal{P}}\left(\boldsymbol{X}\left(t_{0}\right)-\boldsymbol{x}_{0}\right)\delta\left(t\right).
\end{align}
Expressing the time derivative of $\boldsymbol{\mathcal{P}}\boldsymbol{X}_{L}$
as 
\begin{align}
\boldsymbol{\mathcal{P}}\dfrac{d}{dt}\boldsymbol{X}_{L}\left(\left(t-t_{0}\right)/\varepsilon\right) & =\boldsymbol{\delta}_{L,\varepsilon}\left(t-t_{0}\right),\,t\geq t_{0},
\end{align}
finally gives 
\begin{align}
\lim_{\varepsilon\rightarrow0}\boldsymbol{\mathcal{P}}\dfrac{d}{dt}\boldsymbol{X}_{L}\left(\left(t-t_{0}\right)/\varepsilon\right) & =\lim_{\varepsilon\rightarrow0}\boldsymbol{\delta}_{L,\varepsilon}\left(t-t_{0}\right)=2\boldsymbol{\mathcal{P}}\left(\boldsymbol{X}\left(t_{0}\right)-\boldsymbol{x}_{0}\right)\delta\left(t-t_{0}\right).
\end{align}
A similar discussion for the right inner equation yields the equivalent
result
\begin{align}
\lim_{\varepsilon\rightarrow0}\boldsymbol{\mathcal{P}}\dfrac{d}{dt}\boldsymbol{X}_{L}\left(\left(t_{1}-t\right)/\varepsilon\right) & =-2\boldsymbol{\mathcal{P}}\left(\boldsymbol{X}\left(t_{1}\right)-\boldsymbol{x}_{1}\right)\delta\left(t_{1}-t\right).
\end{align}
Finally, the exact solution for the control signal for $\varepsilon=0$
reads as 
\begin{align}
\boldsymbol{u}\left(t\right) & =\begin{cases}
\boldsymbol{\mathcal{B}}^{g}\left(\boldsymbol{x}_{0}\right)\left(\boldsymbol{\dot{X}}\left(t_{0}\right)-\boldsymbol{R}\left(\boldsymbol{x}_{0}\right)\right)+2\boldsymbol{\mathcal{B}}^{g}\left(\boldsymbol{x}_{0}\right)\left(\boldsymbol{X}\left(t_{0}\right)-\boldsymbol{x}_{0}\right)\delta\left(t-t_{0}\right), & t=t_{0},\\
\boldsymbol{\mathcal{B}}^{g}\left(\boldsymbol{X}\left(t\right)\right)\left(\boldsymbol{\dot{X}}\left(t\right)-\boldsymbol{R}\left(\boldsymbol{X}\left(t\right)\right)\right), & t_{0}<t<t_{1},\\
\boldsymbol{\mathcal{B}}^{g}\left(\boldsymbol{x}_{1}\right)\left(\boldsymbol{\dot{X}}\left(t_{1}\right)-\boldsymbol{R}\left(\boldsymbol{x}_{1}\right)\right)-2\boldsymbol{\mathcal{B}}^{g}\left(\boldsymbol{x}_{1}\right)\left(\boldsymbol{X}\left(t_{1}\right)-\boldsymbol{x}_{1}\right)\delta\left(t_{1}-t\right) & t=t_{1}.
\end{cases}
\end{align}
In conclusion, the control diverges at the initial and terminal time,
$t=t_{0}$ and $t=t_{1}$, respectively. The divergence is in form
of a Dirac delta function. The delta kick has a direction in state
space parallel to the jump of the discontinuous state components.
The strength of the delta kick is twice the height of the jump. Inside
the time domain, the control signal is finite and continuous and entirely
given in terms of the outer solution $\boldsymbol{X}\left(t\right)$.

\section*{\label{sec:AppendixD}Appendix E: Analytical results for two-dimensional
dynamical systems}

We state the perturbative solution for optimal trajectory tracking
for the two-dimensional dynamical system 
\begin{align}
\dot{x}\left(t\right) & =a_{0}+a_{1}x\left(t\right)+a_{2}y\left(t\right), & \dot{y}\left(t\right) & =R\left(x\left(t\right),y\left(t\right)\right)+b\left(x\left(t\right),y\left(t\right)\right)u\left(t\right)\label{eq:xyState2}
\end{align}
with state $\boldsymbol{x}=\left(x,y\right)^{T}$, co-state $\boldsymbol{\lambda}=\left(\lambda_{x},\lambda_{y}\right)^{T}$,
and control $u$. The mechanical model from Eq. (\ref{eq:xyState})
of the main text is a special case with $a_{0}=a_{1}=0$ and $a_{2}=1$.
The task is to minimize
\begin{align}
\mathcal{J}\left[\boldsymbol{x}\left(t\right),u\left(t\right)\right] & =\frac{1}{2}\intop_{t_{0}}^{t_{1}}\left(s_{1}\left(x\left(t\right)-x_{d}\left(t\right)\right)^{2}+s_{2}\left(y\left(t\right)-y_{d}\left(t\right)\right)^{2}\right)dt+\frac{\varepsilon^{2}}{2}\intop_{t_{0}}^{t_{1}}dt\left(u\left(t\right)\right)^{2}\label{eq:Functional}
\end{align}
with positive weighting coefficients $s_{1}$ and $s_{2}$. Equation
(\ref{eq:Functional}) corresponds to a choice
\begin{align}
\boldsymbol{\mathcal{S}} & =\left(\begin{array}{cc}
s_{1} & 0\\
0 & s_{2}
\end{array}\right)
\end{align}
for the matrix $\boldsymbol{\mathcal{S}}$ of weighting coefficients.
The necessary optimality conditions are (see Eqs. (\ref{eq:StationarityCondition})-(\ref{eq:AdjointEq})
of the main text) 
\begin{align}
0 & =\varepsilon^{2}u\left(t\right)+b\left(x\left(t\right),y\left(t\right)\right)\lambda_{y}\left(t\right),\label{eq:PontryaginMaximumCondition}\\
\left(\begin{array}{c}
\dot{x}\left(t\right)\\
\dot{y}\left(t\right)
\end{array}\right) & =\left(\begin{array}{c}
y\left(t\right)\\
R\left(x\left(t\right),y\left(t\right)\right)
\end{array}\right)+b\left(x\left(t\right),y\left(t\right)\right)\left(\begin{array}{c}
0\\
u\left(t\right)
\end{array}\right),\label{eq:ControlledTwoDimDynamicalSystem}\\
-\left(\begin{array}{c}
\dot{\lambda}_{x}\left(t\right)\\
\dot{\lambda}_{y}\left(t\right)
\end{array}\right) & =\left(\begin{array}{cc}
a_{1} & \partial_{x}R\left(x\left(t\right),y\left(t\right)\right)+\partial_{x}b\left(x\left(t\right),y\left(t\right)\right)u\left(t\right)\\
a_{2} & \partial_{y}R\left(x\left(t\right),y\left(t\right)\right)+\partial_{y}b\left(x\left(t\right),y\left(t\right)\right)u\left(t\right)
\end{array}\right)\left(\begin{array}{c}
\lambda_{x}\left(t\right)\\
\lambda_{y}\left(t\right)
\end{array}\right)+\left(\begin{array}{c}
s_{1}\left(x\left(t\right)-x_{d}\left(t\right)\right)\\
s_{2}\left(y\left(t\right)-y_{d}\left(t\right)\right)
\end{array}\right),\label{eq:CostateEquations}
\end{align}
which are to be solved with the initial and terminal conditions
\begin{align}
x\left(t_{0}\right) & =x_{0}, & y\left(t_{0}\right) & =y_{0}, & x\left(t_{1}\right) & =x_{1}, & y\left(t_{1}\right) & =y_{1}.
\end{align}
Rearranging the necessary optimality conditions eliminates $\lambda_{y}$
and yields \cite{loeber2015optimal}
\begin{align}
\dot{x}\left(t\right) & =a_{0}+a_{1}x\left(t\right)+a_{2}y\left(t\right),\label{eq:XTransformed}\\
-\dot{\lambda}_{x}\left(t\right) & =a_{1}\lambda_{x}\left(t\right)+s_{1}\left(x\left(t\right)-x_{d}\left(t\right)\right)\nonumber \\
 & -\frac{\varepsilon^{2}}{b\left(x\left(t\right),y\left(t\right)\right)^{2}}\left(\partial_{x}R\left(x\left(t\right),y\left(t\right)\right)+\partial_{x}b\left(x\left(t\right),y\left(t\right)\right)u\left(t\right)\right)\left(\dot{y}\left(t\right)-R\left(x\left(t\right),y\left(t\right)\right)\right),\label{eq:LambdaXTransformed}\\
\varepsilon^{2}\ddot{y}\left(t\right) & =\varepsilon^{2}\dot{x}\left(t\right)\partial_{x}R\left(x\left(t\right),y\left(t\right)\right)+\varepsilon^{2}\dot{y}\left(t\right)\partial_{y}R\left(x\left(t\right),y\left(t\right)\right)-2\varepsilon^{2}\frac{\left(R\left(x\left(t\right),y\left(t\right)\right)-\dot{y}\left(t\right)\right)}{b\left(x\left(t\right),y\left(t\right)\right)}w_{1}\left(x\left(t\right),y\left(t\right)\right)\nonumber \\
 & +b\left(x\left(t\right),y\left(t\right)\right)^{2}\left(a_{2}\lambda_{x}\left(t\right)+s_{2}\left(y\left(t\right)-y_{d}\left(t\right)\right)\right)-\varepsilon^{2}w_{2}\left(x\left(t\right),y\left(t\right)\right)\left(\dot{y}\left(t\right)-R\left(x\left(t\right),y\left(t\right)\right)\right).\label{eq:YTransformed}
\end{align}
 Here, $w_{1}$ and $w_{2}$ denote the abbreviations
\begin{align}
w_{1}\left(x\left(t\right),y\left(t\right)\right) & =\dot{x}\left(t\right)\partial_{x}b\left(x\left(t\right),y\left(t\right)\right)+\dot{y}\left(t\right)\partial_{y}b\left(x\left(t\right),y\left(t\right)\right),\\
w_{2}\left(x\left(t\right),y\left(t\right)\right) & =\partial_{y}R\left(x\left(t\right),y\left(t\right)\right)+\frac{\partial_{y}b\left(x\left(t\right),y\left(t\right)\right)}{b\left(x\left(t\right),y\left(t\right)\right)}\left(\dot{y}\left(t\right)-R\left(x\left(t\right),y\left(t\right)\right)\right).
\end{align}
The outer equations are defined for the ordinary time scale $t$.
The outer solutions are denoted with upper-case letters 
\begin{align}
X\left(t\right) & =x\left(t\right), & Y\left(t\right) & =y\left(t\right), & \Lambda\left(t\right) & =\lambda_{x}\left(t\right).
\end{align}
Expanding Eqs. (\ref{eq:LambdaXTransformed})-(\ref{eq:XTransformed})
up to leading order in $\varepsilon$ yields two linear differential
equations of first order and an algebraic equation, 
\begin{align}
\dot{\Lambda}\left(t\right) & =-a_{1}\Lambda\left(t\right)+s_{1}\left(x_{d}\left(t\right)-X\left(t\right)\right), & \Lambda\left(t\right) & =\frac{s_{2}}{a_{2}}\left(y_{d}\left(t\right)-Y\left(t\right)\right), & \dot{X}\left(t\right) & =a_{0}+a_{1}X\left(t\right)+a_{2}Y\left(t\right).
\end{align}
The solutions for $X$ and $Y$ can be expressed in terms of the state
transition matrix $\boldsymbol{\Phi}\left(t,t_{0}\right)$,
\begin{align}
\left(\begin{array}{c}
X\left(t\right)\\
Y\left(t\right)
\end{array}\right) & =\boldsymbol{\Phi}\left(t,t_{0}\right)\left(\begin{array}{c}
x_{\text{init}}\\
y_{\text{init}}
\end{array}\right)+\intop_{t_{0}}^{t}d\tau\boldsymbol{\Phi}\left(t,\tau\right)\boldsymbol{f}\left(\tau\right),
\end{align}
with
\begin{align}
\boldsymbol{\Phi}\left(t,t_{0}\right) & =\left(\begin{array}{cc}
\cosh\left(\left(t-t_{0}\right)\varphi_{1}\right)+\frac{a_{1}}{\varphi_{1}}\sinh\left(\left(t-t_{0}\right)\varphi_{1}\right) & \frac{a_{2}}{\varphi_{1}}\sinh\left(\left(t-t_{0}\right)\varphi_{1}\right)\\
\frac{a_{2}s_{1}}{s_{2}\varphi_{1}}\sinh\left(\left(t-t_{0}\right)\varphi_{1}\right) & \cosh\left(\left(t-t_{0}\right)\varphi_{1}\right)-\frac{a_{1}}{\varphi_{1}}\sinh\left(\left(t-t_{0}\right)\varphi_{1}\right)
\end{array}\right),
\end{align}
and $\varphi_{1}=\frac{\sqrt{a_{1}^{2}s_{2}+a_{2}^{2}s_{1}}}{\sqrt{s_{2}}}$,
and inhomogeneity
\begin{align}
\boldsymbol{f}\left(t\right) & =\left(\begin{array}{c}
a_{0}\\
\dot{y}_{d}\left(t\right)+a_{1}y_{d}\left(t\right)-\frac{a_{2}s_{1}}{s_{2}}x_{d}\left(t\right)
\end{array}\right).
\end{align}
The constants $x_{\text{init}}$ and $y_{\text{init}}$ must be determined
by matching. Boundary layers occur at both ends of the time domain.
The initial boundary layer at the left end of the time domain is resolved
using the time scale $\tau_{L}=\left(t-t_{0}\right)/\varepsilon$
and scaled solutions
\begin{align}
X_{L}\left(\tau_{L}\right) & =x\left(t\right)=x\left(t_{0}+\varepsilon\tau_{L}\right), & Y_{L}\left(\tau_{L}\right) & =y\left(t\right)=y\left(t_{0}+\varepsilon\tau_{L}\right), & \Lambda_{L}\left(\tau_{L}\right) & =\lambda_{x}\left(t\right)=\lambda_{x}\left(t_{0}+\varepsilon\tau_{L}\right).
\end{align}
The inner left equations are
\begin{align}
Y_{L}''\left(\tau_{L}\right) & =Y_{L}'\left(\tau_{L}\right)^{2}\frac{\partial_{y}b\left(x_{0},Y_{L}\left(\tau_{L}\right)\right)}{b\left(x_{0},Y_{L}\left(\tau_{L}\right)\right)}+s_{2}b\left(x_{0},Y_{L}\left(\tau_{L}\right)\right){}^{2}\left(Y_{L}\left(\tau_{L}\right)-y_{\text{init}}\right), & \Lambda_{L}'\left(\tau_{L}\right) & =0, & X_{L}'\left(\tau_{L}\right) & =0.\label{eq:InnerLeft}
\end{align}
The differential equations do not involve the nonlinearity $R$. As
long as $b$ depends on $Y_{L}$, this is a nonlinear equation. Because
it is autonomous, it can be transformed to a first order ODE \cite{loeber2015optimal}
\begin{align}
Y_{L}'\left(\tau_{L}\right) & =\sqrt{s_{2}}\left(y_{\text{init}}-Y_{L}\left(\tau_{L}\right)\right)\left|b\left(x_{0},Y_{L}\left(\tau_{L}\right)\right)\right|, & Y_{L}\left(0\right) & =y_{0}.\label{eq:EqYL}
\end{align}
An analytical solution of Eq. (\ref{eq:EqYL}) for arbitrary functions
$b$ does not exist in closed form. If $b\left(x,y\right)=b\left(x\right)$
does not depend on $y$, Eq. (\ref{eq:EqYL}) is linear and has the
solution
\begin{align}
Y_{L}\left(\tau_{L}\right) & =y_{\text{init}}+\exp\left(-\sqrt{s_{2}}\left|b\left(x_{0}\right)\right|\tau_{L}\right)\left(y_{0}-y_{\text{init}}\right).
\end{align}
A treatment analogous to the left boundary layer is performed to resolve
the boundary layer at the right end of the time domain. The relevant
time scale is $\tau_{R}=\left(t_{1}-t\right)/\varepsilon$, and the
scaled solutions are defined as
\begin{align}
X_{R}\left(\tau_{R}\right) & =x\left(t\right)=x\left(t_{1}-\varepsilon\tau_{R}\right), & Y_{R}\left(\tau_{R}\right) & =y\left(t\right)=y\left(t_{1}-\varepsilon\tau_{R}\right), & \Lambda_{R}\left(\tau_{R}\right) & =\lambda_{x}\left(t\right)=\lambda_{x}\left(t_{1}-\varepsilon\tau_{R}\right).
\end{align}
The inner right equations are
\begin{align}
Y_{R}''\left(\tau_{R}\right) & =Y_{R}'\left(\tau_{R}\right)^{2}\frac{\partial_{y}b\left(x_{1},Y_{R}\left(\tau_{R}\right)\right)}{b\left(x_{1},Y_{R}\left(\tau_{R}\right)\right)}+b\left(x_{1},Y_{R}\left(\tau_{R}\right)\right){}^{2}s_{2}\left(Y_{R}\left(\tau_{R}\right)-y_{\text{end}}\right), & \Lambda_{R}'\left(\tau_{R}\right) & =0, & X_{R}'\left(\tau_{R}\right) & =0.
\end{align}
These equations are identical in form to the left inner equations
Eqs. (\ref{eq:InnerLeft}). With the analogous considerations as for
the left inner equations, see Eq. (\ref{eq:EqYL}), the solution to
$Y_{R}\left(\tau_{R}\right)$ is given by the first order ODE 
\begin{align}
Y_{R}'\left(\tau_{R}\right) & =\sqrt{s_{2}}\left(y_{\text{end}}-Y_{R}\left(\tau_{R}\right)\right)\left|b\left(x_{1},Y_{R}\left(\tau_{R}\right)\right)\right|, & Y_{R}\left(0\right) & =y_{1}.\label{eq:EqYR}
\end{align}
Finally, the matching constants $x_{\text{init}}$, and $y_{\text{end}}$
are given by $x_{\text{init}}=x_{0}$ and $y_{\text{end}}=Y\left(t_{1}\right)$.
The last matching constant $y_{\text{init}}$ is determined from $X\left(t_{1}\right)=x_{1}$
and yields 
\begin{align}
y_{\text{init}} & =\text{csch}\left(\varphi_{1}\left(t_{1}-t_{0}\right)\right)\left(\intop_{t_{0}}^{t_{1}}\left(\frac{a_{2}s_{1}}{s_{2}}x_{d}\left(\tau\right)-a_{1}y_{d}\left(\tau\right)\right)\sinh\left(\varphi_{1}\left(t_{1}-\tau\right)\right)d\tau-\varphi_{1}\intop_{t_{0}}^{t_{1}}y_{d}\left(\tau\right)\cosh\left(\varphi_{1}\left(t_{1}-\tau\right)\right)d\tau\right)\nonumber \\
 & +\frac{1}{a_{2}}\left(a_{2}y_{d}\left(t_{0}\right)-a_{1}x_{0}-a_{0}\right)-\frac{\text{csch}\left(\varphi_{1}\left(t_{1}-t_{0}\right)\right)}{a_{2}\varphi_{1}}\left(\cosh\left(\varphi_{1}\left(t_{1}-t_{0}\right)\right)\left(a_{0}a_{1}+\varphi_{1}^{2}x_{0}\right)-a_{0}a_{1}-\varphi_{1}^{2}x_{1}\right).
\end{align}
The composite solution is
\begin{align}
x_{\text{comp}}\left(t\right) & =X\left(t\right),\label{eq:xFinalSol}\\
y_{\text{comp}}\left(t\right) & =Y\left(t\right)+Y_{L}\left(\left(t-t_{0}\right)/\varepsilon\right)-y_{\text{init}}+Y_{R}\left(\left(t_{1}-t\right)/\varepsilon\right)-y_{\text{end}},\label{eq:yFinalSol}\\
\lambda_{\text{comp}}\left(t\right) & =\Lambda\left(t\right).\label{eq:lambdaFinalSol}
\end{align}
The control solution is given in terms of the composite solution as
\begin{align}
u\left(t\right) & =\frac{1}{b\left(x_{\text{comp}}\left(t\right),y_{\text{comp}}\left(t\right)\right)}\left(\dot{y}_{\text{comp}}\left(t\right)-R\left(x_{\text{comp}}\left(t\right),y_{\text{comp}}\left(t\right)\right)\right).\label{eq:controlFinalSol}
\end{align}
\bibliographystyle{apsrev4-1}
\bibliography{literature}

\end{document}